\documentclass{siamltex}

\usepackage{amsmath,amsfonts,graphicx,fullpage,pdflscape,color,multirow}
\usepackage{subfig}
\usepackage[pdfborder={0 0 0 0}]{hyperref}

\newcommand{\clawpack}{{\sc clawpack}}
\newcommand{\amrclaw}{{\sc amrclaw}}

\newcommand{\archivelink}{\url{https://bitbucket.org/grady\_lemoine/poro-2d-cartesian-archive}}

\newcommand{\beq}{\begin{eqnarray}}
\newcommand{\eeq}{\end{eqnarray}}
\newcommand{\beqs}{\begin{eqnarray*}}
\newcommand{\eeqs}{\end{eqnarray*}}
\newcommand{\bC}{\mathbf{C}}
\newcommand{\bS}{\mathbf{S}}
\newcommand{\btau}{\mbox{\boldmath $\tau$}}
\newcommand{\bA}{\mbox{$\mathbf{A}$}}
\newcommand{\bB}{\mbox{$\mathbf{B}$}}
\newcommand{\bD}{\mbox{$\mathbf{D}$}}
\newcommand{\bT}{\mbox{$\mathbf{T}$}}
\newcommand{\bG}{\mbox{$\mathbf{G}$}}
\newcommand{\bGamma}{\mbox{$\mathbf{\Gamma}$}}

\newcommand{\bU}{\mbox{$\mathbf{U}$}}
\newcommand{\bu}{\mbox{$\mathbf{u}$}}
\newcommand{\bv}{\mbox{$\mathbf{v}$}}
\newcommand{\bw}{\mbox{$\mathbf{w}$}}
\newcommand{\bq}{\mbox{$\mathbf{q}$}}

\newcommand{\balpha}{\mbox{\boldmath $\alpha$}}
\newcommand{\bbeta}{\mbox{\boldmath $\beta$}}
\newcommand{\bP}{\mbox{$\mathbf{P}$}}
\newcommand{\bQ}{\mbox{$\mathbf{Q}$}}
\newcommand{\bR}{\mbox{$\mathbf{R}$}}
\newcommand{\bK}{\mbox{$\mathbf{K}$}}
\newcommand{\be}{\mbox{$\mathbf{e}$}}
\newcommand{\ep}{\epsilon}
\newcommand{\pr}{\partial}
\newcommand{\B}{\mathbf}
\newcommand{\warn}[1]{}

\author{Grady\ I.\ Lemoine\footnotemark[1]\ \footnotemark[3] \and M.\
Yvonne\ Ou\footnotemark[2] \and Randall\ J.\ LeVeque\footnotemark[1]}

\begin{document}

\title{High-Resolution Finite Volume Modeling of Wave Propagation in
Orthotropic Poroelastic Media}
\maketitle

\renewcommand{\thefootnote}{\fnsymbol{footnote}}

\footnotetext[1]{Department of Applied Mathematics, University of
Washington, Guggenheim Hall Box 352420, Seattle, WA 98195.
Supported in part by NIH grant
5R01AR53652-2 and NSF grant DMS-0914942}
\footnotetext[2]{Department of Mathematical Sciences, University of
Delaware, 501 Ewing Hall, Newark, DE 19716.
Supported in part by NSF-DMS Mathematical Biology Grant 0920852 and NIH CBER PILOT Grant 322 159}
\footnotetext[3]{Corresponding author, email: {\tt gl@uw.edu}}

\begin{abstract}
Poroelasticity theory models the dynamics of porous, fluid-saturated
media.  It was pioneered by Maurice Biot in the 1930s through 1960s,
and has applications in several fields, including geophysics and
modeling of {\em in vivo} bone.  A wide variety of methods have been
used to model poroelasticity, including finite difference, finite
element, pseudospectral, and discontinuous Galerkin methods.  In this
work we use a Cartesian-grid high-resolution finite volume method to
numerically solve Biot's equations in the time domain for orthotropic
materials, with the stiff relaxation source term in the equations
incorporated using operator splitting.  This class of finite volume
method has several useful properties, including the ability to use
wave limiters to reduce numerical artifacts in the solution, ease of
incorporating material inhomogeneities, low memory overhead, and an
explicit time-stepping approach.  To the authors' knowledge, this is
the first use of high-resolution finite volume methods to model
poroelasticity.  The solution code uses the \clawpack{} finite volume
method software, which also includes block-structured adaptive mesh
refinement in its \amrclaw{} variant.  We present convergence results
for known analytic plane wave solutions, achieving second-order
convergence rates outside of the stiff regime of the system.  Our
convergence rates are degraded in the stiff regime, but we still
achieve similar levels of error on the finest grids examined.  We also
demonstrate good agreement against other numerical results from the
literature.  To aid in reproducibility, we provide all of the
code used to generate the results of this paper, at \archivelink{}.
\end{abstract}


\section{Introduction}
\label{sec:intro}

Poroelasticity theory is a homogenized model for solid porous media
containing fluids that can flow through the pore structure.  This
field was pioneered by Maurice A. Biot, who developed his theory of
poroelasticity from the 1930s through the 1960s; a summary of much of
Biot's work can be found in his 1956 and 1962 papers~\cite{biot:56-1,
  biot:56-2, biot:62}.  Biot theory uses linear elasticity to describe
the solid portion of the medium (often termed the {\em skeleton} or
{\em matrix}), linearized compressible fluid dynamics to describe the
fluid portion, and Darcy's law to model the aggregate motion of the
fluid through the matrix.  While it was originally developed to model
fluid-saturated rock and soil, Biot theory has also been used in
underwater acoustics~\cite{bgwx:2004, gilbert-lin:1997,
  gilbert-ou:2003}, and to describe wave propagation in {\em in vivo}
bone~\cite{cowin:bone-poro, cowin-cardoso:fabric-2010,
  gilbert-guyenne-ou:bone-2012}.

Biot theory predicts rich and complex wave phenomena within
poroelastic meterials.  Three different types of waves appear: fast P
waves analogous to standard elastic P waves, in which the fluid and
matrix show little relative motion, and typically
compress or expand in phase with each other; shear waves analogous to
elastic S waves; and slow P waves, where the fluid expands while the
solid contracts, or vice versa.  The slow P waves exhibit substantial
relative motion between the solid and fluid compared to waves of the
other two types.  The viscosity of the fluid dissipates poroelastic waves as
they propagate through the medium, with the fast P and S waves being
lightly damped and the slow P wave strongly damped.  The viscous
dissipation also causes slight dispersion in the fast P and S waves,
and strong dispersion in the slow P wave.

A variety of different numerical approaches have been used to model
poroelasticity.  Carcione, Morency, and Santos provide a thorough
review of the previous
literature~\cite{carcione-morency-santos:comp-poro-rev}.  The earliest
numerical work in poroelasticity seems to be that of
Garg~\cite{garg-nayfeh-good:poro-1974}, using a finite difference
method in 1D.  Finite difference and pseudospectral methods have
continued to be popular since then, with further work by
Mikhailenko~\cite{mikhailenko:poro-1985},
Hassanzadeh~\cite{hassanzadeh:poro-1991}, Dai et
al.~\cite{dai-vafidis-kanasewich:poro-vel-stress}, and more recently
Chiavassa and Lombard~\cite{chiavassa-lombard:poro-fdm-2011}, among
others.  Finite element approaches began being used in the 1980s, with
Santos and Ore\~na's work~\cite{santos-orena:poro-fem} being one of
the first.  Boundary element methods have also been used, such as in
the work of Attenborough, Berry, and
Chen~\cite{attenborough-berry-chen:scattering}.  Spectral element
methods have also been used in both the frequency
domain~\cite{degrande-deroeck:poro-spectral-freq} and the time
domain~\cite{morency-tromp:poro-spectral-time}.  With the recent rise
of discontinuous Galerkin methods, DG has been applied to
poroelasticity in several works, such as that of de la Puente et
al.~\cite{delapuente-dumbser-kaser-igel:poro-dg}.  There have also
been semi-analytical approaches to solving the poroelasticity
equations, such as that of Detournay and
Cheng~\cite{detournay-cheng:borehole}, who analytically obtain a
solution in the Laplace transform domain, but are forced to use an
approximate inversion procedure to return to the time domain.
Finally, there has been significant work on inverse problems in
poroelasticity, for which various forward solvers have been used; of
particular note is the paper of Buchanan, Gilbert, and
Khashanah~\cite{buchanan-gilbert-khashanah:bone-params-lofreq-2004},
who used the finite element method (specifically the FEMLAB software
package) to obtain time-harmonic solutions for cancellous bone as part
of an inversion scheme to estimate poroelastic material parameters,
and the later papers of Buchanan and
Gilbert~\cite{buchanan-gilbert:bone-params-hifreq-1-2007,
  buchanan-gilbert:bone-params-hifreq-2-2007}, where the authors
instead used numerical contour integration of the Green's function.
Numerical work in the 1970s and 1980s focused on isotropic
poroelasticity, with the earliest work on anisotropic poroelasticity
being by Carcione in 1996~\cite{carcione:poro-aniso-1996}.

A major theme in time-domain numerical modeling of poroelasticity has
been the difficulty of handling the viscous dissipation term, which
has its own intrinsic time scale and causes the poroelasticity system
to be stiff, at least if low-frequency waves are being considered.
(The time scales associated with dissipation are independent of
frequency, so at higher frequencies there is less separation between
them and the time scale associated with wave motion --- in other
words, the system is not stiff if sufficiently high-frequency waves
are being considered.)  This viscous dissipation is particularly
problematic for the slow P wave.  While it can still be addressed if
Biot's equations are solved as a unified system (e.g.\ with an
implicit time-integration method), since the viscous dissipation term
is easy to solve analytically in isolation, operator splitting
approaches have also been a popular way to address this issue; we
adopt this approach as well.  Carcione and Quiroga-Goode used operator
splitting in conjunction with a pseudospectral
method~\cite{carcione-quirogagoode:biot-fdm-opsplit}, and Chiavassa
and Lombard use it with a finite difference
approach~\cite{chiavassa-lombard:poro-fdm-2011}.  De la Puente et
al.\ also investigate an operator splitting approach, and encounter
difficulties in obtaining a fast rate of convergence due to the
stiffness of the dissipation term, which pushed the problems they
investigated toward the diffusive limit of the Biot
system~\cite{delapuente-dumbser-kaser-igel:poro-dg}.

Our work in this paper solves a velocity-stress formulation of Biot
linear orthotropic poroelasticity theory using Cartesian-grid
high-resolution finite volume methods.  These methods are
memory-efficient explicit techniques designed to model hyperbolic
systems, and can include wave limiters designed to reduce the effect
of numerical artifacts in the solution.  Since they are based on
solving Riemann problems at grid interfaces, it is also
straightforward to include material inhomogeneities between cells
using these methods.  To our knowledge, this is the first use of
finite volume methods to model poroelasticity.  We employ the
\clawpack{} finite volume method package~\cite{claw.org}, which offers
the use of operator splitting to model viscous dissipation, and
includes optional Berger-Colella-Oliger block-structured adaptive mesh
refinement~\cite{mjb-rjl:amrclaw} to improve solution efficiency for
larger problems.  We verify our code both against known analytical
solutions and against other numerical solutions from the
poroelasticity literature.

\section{Poroelasticity theory for transversely isotropic materials}
\label{sec:poroelasticity}

Biot's equations of poroelasticity are a complicated system of PDEs
that exhibit rich and varied behaviors.  The reader is encouraged to
refer to a detailed treatment of the subject, such as chapter 7 of
Carcione's book~\cite{carcione:wave-book}, but we also provide an
overview here.  After a review of the basic relations modeling the
behavior of a poroelastic medium in sections \ref{sec:stress-strain}
through \ref{sec:eq-plane-strain}, section
\ref{sec:first-order-system} presents the linear first-order system of
PDEs that forms the basis of our numerical model.  Section
\ref{sec:energy-norm} then defines an
energy norm for our state vector that solves some of the scaling
issues associated with modeling poroelasticity in SI units.  The matrix
created for this energy norm allows us to explore some useful
properties of the model: section \ref{sec:hyperbolicity} provides a
concise proof that our first-order system is hyperbolic, and section
\ref{sec:entropy-function} exhibits a strictly convex entropy function
for the system, which has implications for the correctness of our
numerical solutions that we explore later in section
\ref{sec:stiff-subcharacteristic}.

\subsection{Stress-strain relation}
\label{sec:stress-strain}

We assume the constituent material of the solid matrix is isotropic
and the anisotropy of the solid matrix results purely from the
microstructure.  We also assume that the anisotropy has a specific
form --- that the medium is orthotropic, possessing three orthogonal
planes of symmetry, and is isotropic with respect to some axis of
symmetry.  This type of anisotropy is common in engineering
composites~\cite{gibson:composites}, and in biological
materials~\cite{cowin-mehrabadi:elastic-symmetry}, as well as being
present in certain types of stone.  Let the $z$-axis be this axis of
symmetry.  The elastic stiffness tensor $\bC$ of such an orthotropic,
transversely isotropic medium contains five independent components.
In its principal axes, and using shorthand notation, $\bC$ can be
arranged as \beq \bC=
\begin{pmatrix}
c_{11} & c_{12} & c_{13} & 0 & 0 & 0\\
c_{12} & c_{11} & c_{13} & 0 & 0 & 0\\
c_{13} & c_{13} & c_{33} & 0 & 0 & 0 \\
0 & 0 & 0 & c_{55} & 0 & 0 \\
0 & 0 & 0 & 0 & c_{55} & 0 \\
0 & 0 & 0 & 0 & 0 &\frac{c_{11}-c_{12}}{2} 
\end{pmatrix},
\eeq
with the stress tensor and engineering strains arranged into vectors
$\btau$ and $\be$ of the form
\beq
\btau=\begin{pmatrix}
       \tau_{11} & \tau_{22} & \tau_{33} & \tau_{23} & \tau_{13} & \tau_{12}
      \end{pmatrix}^T
\quad \text{and} \quad
\be=\begin{pmatrix}
       \ep_{11} & \ep_{22} & \ep_{33} & 2\ep_{23} & 2\ep_{13} & 2\ep_{12}
      \end{pmatrix}^T.
\eeq
Note the factor of 2 applied to the shear strains to convert from the
tensor strain $\ep_{ij} = \frac{1}{2}(\pr_i u_j + \pr_j u_i)$ to the
engineering strain. With this shorthand notation, the
stress-strain relation is
\begin{equation}
 \btau=\bC \be .
\end{equation}

\subsection{Energy densities and the dissipation potential}

\label{sec:energy}

One useful property of the poroelasticity system is that it admits an
energy density, from which we can define an energy norm that we will
use extensively.  

This section is mainly based on Biot's 1956 papers~\cite{biot:56-1,
  biot:56-2} and chapter 7 of Carcione's
book~\cite{carcione:wave-book}. All formulations are in terms of the
following variables:
\begin{enumerate}
\item
$\bu$, the displacement vector of the solid matrix
\item
$\bw := \phi(\bU-\bu)$, the relative motion of the fluid scaled by the
  porosity, where $\bU$ is the displacement vector of the pore fluid
  and $\phi$ is the porosity of the medium
\item
$\zeta := -\nabla\cdot\bw$, the variation in fluid content
\end{enumerate}

\subsubsection{Strain energy}
In terms of the undrained elasticity tensor $\bC^{u}$ and the strain
components of the solid matrix $e_{ii}^{(m)}:=\pr_i u_i$ and $e_{ij}^{(m)}:=\pr_i
u_j +\pr_j u_i$, $i\ne j$, the strain energy density of the Biot
model for 3D transversely isotropic materials, with $z$-axis being the
axis of symmetry, is given by
\begin{multline} \label{eq:strain-energy}
2V = c_{11}^u\left({e^{(m)}_{11}}^2+{e^{(m)}_{22}}^2\right)+c_{33}^u {e^{(m)}_{33}}^2+2 c_{12}^u e^{(m)}_{11}e^{(m)}_{22}+2c_{13}^u\left(e^{(m)}_{11}+e^{(m)}_{22}\right) e^{(m)}_{33}\\
+c_{55}^u({e^{(m)}_{23}}^2+{e^{(m)}_{13}}^2)+c_{66}^u
{e^{(m)}_{12}}^2-2 \alpha_1 M \left(e^{(m)}_{11}+e^{(m)}_{22}\right)
\zeta-2\alpha_3 M e^{(m)}_{33} \zeta +M \zeta^2 ,
\end{multline}
where $c_{66}^u=\frac{c_{11}^u-c_{12}^u}{2}$ for a transversely
isotropic material and the undrained elastic constants $c_{ij}^u$ are
related to those of the dry matrix, $c_{ij}^{(m)}$, via
\begin{align}
c^u_{ij}&=c_{ij}^{(m)}+M\alpha_i\alpha_j, \quad i,j=1,\dotsc,6\\
\balpha&:=(\alpha_1, \alpha_1, \alpha_3, 0,0,0)\\
\alpha_1&:=1-\frac{c_{11}^{(m)}+c_{12}^{(m)}+c_{13}^{(m)}}{3K_s}\\
\alpha_3&:=1-\frac{2c_{13}^{(m)}+c_{33}^{(m)}}{3K_s}\\
M&:=\frac{K_s^2}{K_s\left[1+\phi(K_s/K_f -1)\right]-\left(2c_{11}^{(m)}+c_{33}^{(m)}+2c_{12}^{(m)}+4c_{13}^{(m)}\right)/9}.
\end{align}
Here $K_s$ and $K_f$ are the bulk moduli of the constituent material
of the solid matrix and of the pore fluid,
respectively. The total stress (solid matrix plus pore pressure)
acting on a volume element of the medium are given by the derivatives
of strain energy with respect to the associated strains,
\begin{equation}
\tau_{ij}=\frac{\pr V}{\pr e_{ij}^{(m)}},
\end{equation}
or, in short notation,
\beq
\tau_{I}=\sum_{J=1}^6  c^u_{IJ}e_J^{(m)}- M \alpha_I \zeta.
\label{eq:short-stress}
\eeq
Similarly, the pore pressure $p$ is
\begin{equation} \label{eq:pressure}
p=\frac{\pr V}{\pr \zeta}= M\left(\zeta-\sum_{j=1}^3 \alpha_j e_{jj}^{(m)}\right).
\end{equation}

Most of the work in this paper is for 2D plane-strain conditions in
the $x$-$z$ plane.  For these conditions, we will later need a
function $\tilde{V}$ that maps from the pressure and the in-plane stress
components $\tau_{11}$, $\tau_{13}$, $\tau_{33}$ at a point to the
strain energy density at that point.  To obtain this, we will first
derive the strains as a function of stress for plane-strain conditions
by setting $e_{12}^{(m)} = e_{22}^{(m)} = e_{23}^{(m)} = 0$ in equations
\eqref{eq:short-stress} and \eqref{eq:pressure}, then solving for the remaining
strains and the variation in fluid content:
\begin{equation} \label{eq:plane-strain-from-stress}
  \be_{\text{plane}} = \bS_{\text{plane}} \btau_{\text{plane}},
\end{equation}
where the matrices $\be_{\text{plane}}$, $\bS_{\text{plane}}$, and
$\btau_{\text{plane}}$ are
\begin{equation} \label{eq:plane-strain-constitutive-vecs}
  \be_{\text{plane}} = \begin{pmatrix} e_{11}^{(m)} & e_{33}^{(m)} & e_{13}^{(m)} & \zeta \end{pmatrix}^T, \quad
  \btau_{\text{plane}} = \begin{pmatrix} \tau_{11} & \tau_{33} & \tau_{13} & p \end{pmatrix}^T,
\end{equation}
and
\begin{equation} \label{eq:plane-strain-constitutive-mat}
  \bS_{\text{plane}} = \begin{pmatrix}
    \frac{c^{(m)}_{33}}{c^{(m)}_{11}c^{(m)}_{33} - \left(c^{(m)}_{13}\right)^2} &
    -\frac{c^{(m)}_{13}}{c^{(m)}_{11}c^{(m)}_{33} - \left(c^{(m)}_{13}\right)^2} & 0 &
    \frac{\alpha_1 c^{(m)}_{33} - \alpha_3 c^{(m)}_{13}}{c^{(m)}_{11}c^{(m)}_{33} - \left(c^{(m)}_{13}\right)^2} \\
    -\frac{c^{(m)}_{13}}{c^{(m)}_{11}c^{(m)}_{33} - \left(c^{(m)}_{13}\right)^2} &
    \frac{c^{(m)}_{11}}{c^{(m)}_{11}c^{(m)}_{33} - \left(c^{(m)}_{13}\right)^2} & 0 &
    \frac{\alpha_3 c^{(m)}_{11} - \alpha_1 c^{(m)}_{13}}{c^{(m)}_{11}c^{(m)}_{33} - \left(c^{(m)}_{13}\right)^2} \\
    0 & 0 & \frac{1}{c^{(m)}_{55}} & 0 \\
    \frac{\alpha_1 c^{(m)}_{33} - \alpha_3 c^{(m)}_{13}}{c^{(m)}_{11}c^{(m)}_{33} - \left(c^{(m)}_{13}\right)^2} &
    \frac{\alpha_3 c^{(m)}_{11} - \alpha_1 c^{(m)}_{13}}{c^{(m)}_{11}c^{(m)}_{33} - \left(c^{(m)}_{13}\right)^2} & 0 &
    \frac{1}{M} + \frac{\alpha_1^2 c^{(m)}_{33} + \alpha_3^2 c^{(m)}_{11} -
      2\alpha_1\alpha_3 c^{(m)}_{13}}{c^{(m)}_{11}c^{(m)}_{33} - \left(c^{(m)}_{13}\right)^2}
  \end{pmatrix}.
\end{equation}
Substituting \eqref{eq:plane-strain-from-stress} through
\eqref{eq:plane-strain-constitutive-mat} into \eqref{eq:strain-energy}
yields the strain energy as a function of stress for plane strain
conditions, which can be expressed as the quadratic form
\begin{equation} \label{eq:vtilde}
  \tilde{V} := \frac{1}{2} \btau_{\text{plane}}^T \bS_{\text{plane}} \btau_{\text{plane}}
\end{equation}
Note the use of the drained elastic coeffcients ($c^{(m)}_{11}$, etc.)
rather than the undrained coefficients ($c_{11}^u$).  We expect this
quadratic form to be positive-definite on physical grounds --- if it
were not, then it would be possible to deform the medium or change its
fluid content without doing work.

\subsubsection{Kinetic energy and the dissipation potential}
For anisotropic poroelastic media, the kinetic energy density has the
form
\begin{equation}
T = \frac{1}{2} \left( \dot{\bu}^T \bP\dot{\bu}+2
\dot{\bu}^T \bR \dot{\bU} +\dot{\bU}^T \bT \dot{\bU}\right),
\end{equation}
where $\bR$ is the
induced mass matrix. Assume all the three matrices can be diagonalized
in the same coordinate system so that $\bP=\diag(a_1,a_2,a_3)$,
$\bR=\diag(r_1,r_2,r_3)$ and $\bT=\diag(t_1,t_2,t_3)$. By looking at
the special case of no relative motion between fluid and solid, it can
be shown that the mass coefficients satisfy the two equations (given
as (7.169) in Carcione~\cite{carcione:wave-book})
\begin{equation}
a_i+r_i=(1-\phi)\rho_s,\quad
r_i+t_i=\phi\rho_f.
\end{equation}
Defining the velocity variables  
\begin{equation}
\bv := \dot{\bu}, \quad \bq := \dot{\bw},
\end{equation}
the kinetic energy density can be expressed as
\begin{equation} \label{eq:kinetic-energy}
\begin{aligned}
T&=\frac{1}{2}\sum_{i=1}^3 \left[ (1-\phi)\rho_s {v_i}^2 -r_i ({v_i} -\dot{U_i})^2+\phi \rho_f \dot{U_i}^2 \right]\\
&=\frac{1}{2}\sum_{i=1}^3 \left[ \rho v_i^2+2\rho_f q_i
  v_i+\left(\frac{\rho_f \phi - r_i}{\phi^2}\right) q_i^2
  \right] =: \tilde{T}(\bv,\bq),
\end{aligned}
\end{equation}
where $\rho_s$ and $\rho_f$ are the constituent solid density and pore
fluid density, respectively and $\rho=(1-\phi)\rho_s+\phi \rho_f$ is
the bulk density of the medium.
The induced mass parameters $r_i$ are related to the tortuosity $T_i$ by 
\begin{equation}
r_i = \phi \rho_f (1-T_i).
\end{equation}

Assuming the pore fluid flow is of the Poiseuille type, the
dissipation potential $\Phi_D$ in an anisotropic medium in terms of the
dynamic viscosity $\eta$ of pore fluid and the permeability tensor
$\bK$ is
\begin{equation}
\Phi_D=\frac{1}{2} \sum_{i,j=1}^3 \eta {\left(\bK^{-1}\right)}_{ij}
q_i q_j .
\end{equation}
Assume that $\bK$ has the same principal directions as $\bP$, $\bR$
and $\bT$ with eigenvalues $\kappa_1,\kappa_2, \kappa_3$.  Then we have
\beq
\Phi_D=\frac{1}{2} \sum_{i=1}^3 \frac{\phi^2\eta}{\kappa_i} (\dot{U_i}-v_i)^2=\frac{1}{2} \sum_{i=1}^3 \frac{\eta}{\kappa_i} q_i^2=: \tilde{\Phi}_D(\bq)
\eeq

\subsection{Equations of motion} \label{sec:eom}
The equations of motion for the solid part are given in terms of the
energy densities and dissipation potential as
\begin{equation}
\pr_t \left( \frac{\pr T}{\pr v_i}\right)+\frac{\pr \Phi_D}{\pr v_i} =
\sum_{j=1}^3\pr_j \left(\frac{\pr V}{\pr e^{(m)}_{ij}} + \phi p \delta_{ij} \right),
\end{equation}
or equivalently
\beq
\pr_t \left( \frac{\pr \tilde{T}}{\pr v_i}-\phi  \frac{\pr
  \tilde{T}}{\pr q_i} \right)-\phi \frac{\pr \tilde{\Phi}_D}{\pr q_i}
= \sum_{j=1}^3\pr_j \left(\frac{\pr V}{\pr e^{(m)}_{ij}} + \phi p \delta_{ij} \right).
\label{eq-motion-solid}
\eeq
Similarly, the equation of motion for the fluid part is
\begin{equation}
\pr_t \left( \frac{\pr T}{\pr \dot{U_i}}\right)+\frac{\pr \Phi_D}{\pr \dot{U_i}} = -\phi \sum_{j=1}^3\pr_j \left(\frac{\pr V}{\pr\zeta}\delta_{ij} \right),
\end{equation}
or equivalently
\beq
\phi \pr_t \left( \frac{\pr \tilde{T}}{\pr q_i}\right)+\phi \frac{\pr \tilde{\Phi}_D}{\pr q_i} = -\phi \sum_{j=1}^3\pr_j \left(\frac{\pr V}{\pr\zeta}\delta_{ij} \right),
\label{eq-motion-fluid-enegry}
\eeq
which reduces to
\beq
-\pr_i p= \rho_f \dot{v_i}+\left(\frac{\rho_f \phi-r_i}{\phi^2}\right)\dot{q_i}+\frac{\eta}{\kappa_i}q_i,\ i=1,2,3.
\label{eq-motion-fluid}
\eeq

The equations of motion for the fluid-solid composite are obtained by
adding \eqref{eq-motion-solid} with \eqref{eq-motion-fluid-enegry}, giving
\beq
\sum_{j=1}^3 \pr_j\tau_{ij} =\pr_t\left(\frac{\pr\tilde{T}}{\pr v_i}\right)=\rho \dot{v_i}+\rho_f \dot{q_i}.
\label{eq-motion-total}
\eeq

\subsection{Governing equations for plane-strain case}
\label{sec:eq-plane-strain}

Since the material is assumed isotropic in the $x$-$y$
plane, we consider the plane strain problem in the $x$-$z$ plane. The
governing equations are obtained by suppressing the $y$-component
(subscript 2) of $\bU$, $\bu$, $\bw$, $\bq$ and those terms of
$e_{ij}^{(m)}$ with $i=2$ or $j=2$ in $V$,
$T$, $\tilde{T}$, $\Phi_D$, and $\tilde{\Phi}_D$.  (While nonzero
out-of-plane stresses do arise in a plane-strain problem, they do not
produce in-plane motion, and can be ignored for purposes of studying the
in-plane dynamics of the medium.)  We obtain two types of governing equation:

\begin{itemize}
\item Stress-strain relations, obtained by differentiating
  \eqref{eq:pressure} and \eqref{eq:short-stress} for $I = 1,3,5$ with
  respect to time,
\begin{align}
\pr_t \tau_{xx}&=c_{11}^u \pr_x v_x+c_{13}^u\pr_z v_z+\alpha_1 M(\pr_x q_x+\pr_z q_z)+\pr_t s_1\label{stress-strain-1}\\
\pr_t \tau_{zz}&=c_{13}^u \pr_x v_x+c_{33}^u\pr_z v_z+\alpha_3 M(\pr_x q_x+\pr_z q_z)+\pr_t s_3\\
\pr_t \tau_{xz}&=c_{55}^u (\pr_z v_x+\pr_x v_z)+\pr_t s_5\\
\pr_t p &= -\alpha_1 M \pr_x v_x -\alpha_3 M \pr_z v_z -M(\pr_x q_x+\pr_z q_z)+\pr_t s_f,\label{stress-strain-4}
\end{align}
where $s_1$, $s_3$, $s_5$, and $s_f$ are the solid and fluid external sources.
\item Equations of motion
\begin{align}
\rho \pr_t v_x + \rho_f \pr_t q_x &= \pr_x \tau_{xx}+\pr_z\tau_{xz}\label{motion-x}\\
\rho \pr_t v_z + \rho_f \pr_t q_z &= \pr_x \tau_{xz}+\pr_z\tau_{zz}\label{motion-z}\\
\rho_f \pr_t v_x + m_1 \pr_t q_x +\left(\frac{\eta}{\kappa_1}\right)q_x &= -\pr_x p\label{darcy-x}\\
\rho_f \pr_t v_z + m_3 \pr_t q_z +\left(\frac{\eta}{\kappa_3}\right)q_z &= -\pr_z p\label{darcy-z},
\end{align}
where $m_i:=\frac{\rho_f \phi -r_i}{\phi^2} = \frac{\rho_fT_i}{\phi}$.
\end{itemize}
\subsection{Governing equations as a linear first-order system}
\label{sec:first-order-system}
Solving \eqref{motion-x} and \eqref{darcy-x} for $\pr_t v_x$ and $\pr_t q_x$, we obtain
\begin{align}
\pr_t v_x &= \frac{1}{\Delta_1}\left(m_1 \pr_x \tau_{xx}+m_1 \pr_z \tau_{xz}+\rho_f \pr_x p +\rho_f\frac{\eta}{\kappa_1}q_x\right)\label{motion-darcy-1}\\
\pr_t q_x &= \frac{1}{\Delta_1}\left(-\rho_f \pr_x \tau_{xx}-\rho_f  \pr_z \tau_{xz}- \rho \pr_x p -\rho\frac{\eta}{\kappa_1}q_x\right),
\end{align}
where
$\Delta_1:= \rho m_1 -\rho_f^2$.
Similarly, \eqref{motion-z} and \eqref{darcy-z} lead to
\begin{align}
\pr_t v_z &= \frac{1}{\Delta_3}\left(m_3 \pr_x \tau_{xz}+m_3 \pr_z \tau_{zz}+\rho_f \pr_z p +\rho_f\frac{\eta}{\kappa_3}q_z\right)\\
\pr_t q_z &= \frac{1}{\Delta_3}\left(-\rho_f \pr_x \tau_{xz}-\rho_f  \pr_z \tau_{zz}- \rho \pr_z p -\rho\frac{\eta}{\kappa_3}q_z\right),\label{motion-darcy-4}
\end{align}
where
$\Delta_3:= \rho m_3 -\rho_f^2$.
Combining the stress-strain relations \eqref{stress-strain-1}-\eqref{stress-strain-4} with equations \eqref{motion-darcy-1}-\eqref{motion-darcy-4}, we obtain the $8\times 8$ linear first-order system
\beq \label{eq:first-order-system}
\pr_t \bQ + \B{A} \pr_x\bQ+\B{B} \pr_z \bQ = \B{D} \bQ + \pr_t \B{s},
\eeq
where
\begin{align}
\bQ&=\begin{pmatrix} \tau_{xx} & \tau_{zz} & \tau_{xz} & v_x & v_z & p & q_x & q_z\end{pmatrix}^T \\
\B{A}&=-\begin{pmatrix} \label{eq:matrix-A}
0 & 0 &0 &c_{11}^u &0 &0 &\alpha_1 M &0\\
0 & 0 &0 &c_{13}^u &0 &0 &\alpha_3 M &0\\
0 & 0 &0 &0 &c_{55}^u &0 &0  &0\\
\frac{m_1}{\Delta_1} &0 &0 &0 &0 &\frac{\rho_f}{\Delta_1} &0 &0 \\
0 & 0 &\frac{m_3}{\Delta_3} &0 &0 &0 &0 &0 \\
0 & 0 &0 &-\alpha_1 M & 0 &0 &-M &0\\
-\frac{\rho_f}{\Delta_1} &0 &0 &0 &0 &-\frac{\rho}{\Delta_1} &0 &0\\
0 & 0 &-\frac{\rho_f}{\Delta_3} &0 &0 &0 &0 &0
\end{pmatrix}\\
\B{B}&=-\begin{pmatrix} \label{eq:matrix-B}
0 & 0 &0 &0 &c_{13}^u &0 &0 &\alpha_1 M\\
0 & 0 &0 &0 &c_{33}^u &0 &0 &\alpha_3 M\\
0 & 0 &0 &c_{55}^u &0 &0  &0 &0\\
0 &0 &\frac{m_1}{\Delta_1} &0 &0 &0 &0 &0 \\
0 &\frac{m_3}{\Delta_3} &0 &0 &0 &\frac{\rho_f}{\Delta_3} &0 &0 \\
0 & 0 &0 &0 &-\alpha_3 M & 0 &0 &-M \\
0 &0 &-\frac{\rho_f}{\Delta_1} &0 &0  &0 &0 &0\\
0 &-\frac{\rho_f}{\Delta_3} &0 &0 &0 &-\frac{\rho}{\Delta_3} &0 &0
\end{pmatrix}\\
\B{D}&=\begin{pmatrix} \label{eq:matrix-D}
0 & 0 &0 &0 &0 &0 &0 &0\\
0 & 0 &0 &0 &0 &0 &0 &0\\
0 & 0 &0 &0 &0 &0 &0 &0\\
0 & 0 &0 &0 &0 &0 &\frac{\rho_f\eta}{\Delta_1\kappa_1} &0 \\
0 & 0 &0 &0 &0 &0 &0 &\frac{\rho_f\eta}{\Delta_3\kappa_3}\\
0 & 0 &0 &0 &0 &0 &0 &0\\
0 & 0 &0 &0 &0 &0 &-\frac{\rho\eta}{\Delta_1\kappa_1} &0\\
0 & 0 &0 &0 &0 &0 &0 &-\frac{\rho\eta}{\Delta_3\kappa_3}
\end{pmatrix},
\end{align}
and 
\begin{equation}
 \B{s}=\begin{pmatrix} s_1 & s_3 & s_5 &0 &0 &s_f &0 &0 \end{pmatrix}^T.
\end{equation}
It is this system \eqref{eq:first-order-system} that forms the basis
for our numerical work.  Note that while the coefficient matrices
$\B{A}$, $\B{B}$, and $\B{D}$ are defined in the material principal
axes, we can extend this system to model media where the principal
axes are different from the global $x$-$z$ axes through an appropriate
transformation of the state variables in $\B{Q}$ that come from vector
and tensor quantities, and application of the chain rule in the
partial derivatives with respect to the spatial variables.  In such
cases, we refer to the principal directions as the $1$ and $3$
axes to distinguish them from the computational $x$ and $z$
axes.

It is worth noting that \eqref{eq:first-order-system} is not just a
generic ``black box'' equation, but one of a very specific type: a
first-order hyperbolic system with a stiff relaxation source term.  We
will prove hyperbolicity in section \ref{sec:hyperbolicity}, and the
source term $\bD \bQ$ shows itself to be of relaxation type by having
only zero and negative eigenvalues --- in the absence of the spatial
derivative terms, it would cause the solution to decay exponentially
toward the null space $\mathcal{N}(\bD)$, and even with the other
terms present we can expect it to keep the solution close to
$\mathcal{N}(\bD)$.  (Whether the relaxation term really is stiff
depends on the other time scales of the particular problem being
solved, but it is stiff for some of the problems
considered here.)  On the subject of time scales, because $\bD$ is
extremely sparse, we can immediately read off the eigenvalues
associated with dissipation in the $1$ and $3$ axes --- respectively,
$-\frac{\rho\eta}{\Delta_1\kappa_1}$ and
$-\frac{\rho\eta}{\Delta_3\kappa_3}$ --- and so define the
characteristic time for decay in each axis as the negative inverse of
these eigenvalues,
\begin{equation} \label{eq:decay-characteristic-times}
  \tau_{d1} := \frac{\Delta_1\kappa_1}{\rho\eta}, \quad
  \tau_{d3} := \frac{\Delta_3\kappa_3}{\rho\eta}.
\end{equation}

\subsection{Energy norm}
\label{sec:energy-norm}

Let $\mathcal{E} := \tilde{T} + \tilde{V}$ be the total mechanical
energy per unit volume in a representative element, where $\tilde{T}$
is the kinetic energy from \eqref{eq:kinetic-energy}, and $\tilde{V}$
is the strain energy from \eqref{eq:vtilde}.  For the
subsequent analysis, we will use its Hessian with respect to the state
variables in $\bQ$, which is the symmetric matrix
\begin{equation} \label{eq:energy-hessian}
  \B{E} = \begin{pmatrix}
    \frac{c^{(m)}_{33}}{c_{11}^{(m)}c_{33}^{(m)} - \left(c_{13}^{(m)}\right)^2} &
    -\frac{c^{(m)}_{13}}{c_{11}^{(m)}c_{33}^{(m)} - \left(c_{13}^{(m)}\right)^2} & 0 & 0 & 0 &
    \frac{\alpha_1 c^{(m)}_{33} - \alpha_3 c^{(m)}_{13}}{c^{(m)}_{11}c^{(m)}_{33} - \left(c_{13}^{(m)}\right)^2} &
    0 & 0 \\
    -\frac{c^{(m)}_{13}}{c^{(m)}_{11}c^{(m)}_{33} - \left(c_{13}^{(m)}\right)^2} &
    \frac{c^{(m)}_{11}}{c^{(m)}_{11}c^{(m)}_{33} - \left(c_{13}^{(m)}\right)^2} & 0 & 0 & 0 &
    \frac{\alpha_3 c^{(m)}_{11} - \alpha_1 c^{(m)}_{13}}{c^{(m)}_{11}c^{(m)}_{33} - \left(c_{13}^{(m)}\right)^2} &
    0 & 0 \\
    0 & 0 & \frac{1}{c^{(m)}_{55}} & 0 & 0 & 0 & 0 & 0 \\
    0 & 0 & 0 & \rho & 0 & 0 & \rho_f & 0 \\
    0 & 0 & 0 & 0 & \rho & 0 & 0 & \rho_f \\
    \frac{\alpha_1 c^{(m)}_{33} - \alpha_3 c^{(m)}_{13}}{c^{(m)}_{11}c^{(m)}_{33} - \left(c_{13}^{(m)}\right)^2} &
    \frac{\alpha_3 c^{(m)}_{11} - \alpha_1 c^{(m)}_{13}}{c^{(m)}_{11}c^{(m)}_{33} - \left(c_{13}^{(m)}\right)^2} &
    0 & 0 & 0 & \frac{1}{M} + \frac{\alpha_1^2 c^{(m)}_{33} + \alpha_3^2 c^{(m)}_{11} -
    2\alpha_1\alpha_3 c^{(m)}_{13}}{c^{(m)}_{11}c^{(m)}_{33} - \left(c_{13}^{(m)}\right)^2} & 0 & 0 \\
    0 & 0 & 0 & \rho_f & 0 & 0 & m_1 & 0 \\
    0 & 0 & 0 & 0 & \rho_f & 0 & 0 & m_3
  \end{pmatrix}.
\end{equation}
We know that $\B{E}$ is a positive-definite matrix because it
is the Hessian of the positive-definite quadratic form $\mathcal{E}$.
In fact, because $\mathcal{E}$ has no linear terms in the state
variables, we can write it compactly in terms of its Hessian as
\begin{equation} \label{eq:energy-from-hessian}
  \mathcal{E} = \frac{1}{2}\bQ^T \B{E} \bQ.
\end{equation}

For many poroelastic materials, the components of $\bQ$ are very badly
scaled relative to each other when expressed in common units -- for
example, waves in the geological materials of Table \ref{tab:matprops}
typically have stress components about seven orders of magnitude
larger than their velocity components when expressed in SI base units.
This makes using the usual vector norms on $\bQ$ problematic, but we can fix
this issue by using $\B{E}$ to define an energy norm,
\begin{equation} \label{eq:energy-norm}
  \|\bQ\|_E := \sqrt{\bQ^H \B{E} \bQ}.
\end{equation}
We define the norm using the Hermitian conjugate-transpose
(superscript $H$) rather than
the simple vector transpose in case we later want to take the energy
norm of a complex vector.  This norm has the advantage of scaling the
elements of $\bQ$ in a physically relevant fashion, and producing a
result that is physically meaningful and has consistent units.
Incidentally, this also lets us write the energy density in the even
more compact form
\begin{equation}
  \mathcal{E} = \frac{1}{2}\|\bQ\|_E^2.
\end{equation}

Note that because $\B{E}$ involves the elastic moduli, density,
etc. of the medium, it is different for different materials in a
heterogeneous domain.  Because the energy norm is derived from the
energy density at a point, however, energy norms computed in different
materials can still be meaningfully compared.

\subsection{Hyperbolicity}
\label{sec:hyperbolicity}

Equipped with the matrix $\B{E}$, we can now prove that the left-hand
side of \eqref{eq:first-order-system} forms a hyperbolic system.

First, consider the system formed by setting the left-hand side of
\eqref{eq:first-order-system} equal to zero,
\begin{equation} \label{eq:first-order-hyperbolic}
  \pr_t \bQ + \bA \pr_x \bQ + \bB \pr_z \bQ = 0.
\end{equation}
This system is hyperbolic if the linear combination $\breve{\bA} = n_x
\bA + n_z \bB$ is diagonalizable with real eigenvalues for all real
$n_x$ and $n_z$.  Suppose $\B{v}$ is an eigenvector of $\breve{\bA}$
with eigenvalue $\lambda$.  Then
\begin{equation} \label{eq:eigenproblem-base}
 \breve{\bA} \B{v} = \lambda \B{v},
\end{equation}
or, multiplying by $\B{E}$,
\begin{equation} \label{eq:eigenproblem-energy}
  \B{E} \breve{\bA} \B{v} = \lambda \B{E} \B{v}.
\end{equation}
Since $\B{E}$ is nonsingular, any pair $(\B{v}, \lambda)$ that
satisfies the generalized eigenproblem \eqref{eq:eigenproblem-energy}
also satisfies the original eigenproblem
\eqref{eq:eigenproblem-base}.

It is not obvious how this transformation is helpful, but if we
examine the component matrices $\B{E} \bA$ and $\B{E} \bB$ of $\B{E}
\breve{\bA}$, after substantial algebra we can discover they are
symmetric:
\begin{equation} \label{eq:matrices-symmetrized}
  \B{E}\bA = \begin{pmatrix}
    0 & 0 & 0 & -1 & 0 & 0 & 0 & 0 \\
    0 & 0 & 0 & 0 & 0 & 0 & 0 & 0 \\
    0 & 0 & 0 & 0 & -1 & 0 & 0 & 0 \\
    -1 & 0 & 0 & 0 & 0 & 0 & 0 & 0 \\
    0 & 0 & -1 & 0 & 0 & 0 & 0 & 0 \\
    0 & 0 & 0 & 0 & 0 & 0 & 1 & 0 \\
    0 & 0 & 0 & 0 & 0 & 1 & 0 & 0 \\
    0 & 0 & 0 & 0 & 0 & 0 & 0 & 0
  \end{pmatrix}, \quad
  \B{E}\bB = \begin{pmatrix}
    0 & 0 & 0 & 0 & 0 & 0 & 0 & 0 \\
    0 & 0 & 0 & 0 & -1 & 0 & 0 & 0 \\
    0 & 0 & 0 & -1 & 0 & 0 & 0 & 0 \\
    0 & 0 & -1 & 0 & 0 & 0 & 0 & 0 \\
    0 & -1 & 0 & 0 & 0 & 0 & 0 & 0 \\
    0 & 0 & 0 & 0 & 0 & 0 & 0 & 1 \\
    0 & 0 & 0 & 0 & 0 & 0 & 0 & 0 \\
    0 & 0 & 0 & 0 & 0 & 1 & 0 & 0
  \end{pmatrix}.
\end{equation}
Thus $\B{E} \breve{\bA}$ is a symmetric matrix, and
\eqref{eq:eigenproblem-energy} is a real symmetric generalized
eigenproblem with a positive-definite matrix on its right-hand side.
As such, it has purely real eigenvalues and a full set of linearly
independent eigenvectors.  $\breve{\bA}$ is therefore diagonalizable
and has pure real eigenvalues, which means that
\eqref{eq:first-order-hyperbolic} is a hyperbolic system.

\subsection{Entropy function} \label{sec:entropy-function}

We can also show that the energy density is a strictly convex entropy
function of the system \eqref{eq:first-order-system}, in a sense
similar to that of Chen, Levermore, and
Liu~\cite{chen-levermore-liu:stiff-relaxation}.  Adapting the
definition of~\cite{chen-levermore-liu:stiff-relaxation} to the
notation used here, a function $\Phi: \mathbb{R}^8 \to \mathbb{R}$ is
a strictly convex entropy function for the system
\eqref{eq:first-order-system} if it satisfies the following
conditions:
\begin{enumerate}
\item $\Phi''(\bQ) (n_x \bA + n_z \bB)$ is symmetric for
  all scalars $n_x$ and $n_z$
\item $(\Phi'(\bQ))^T \B{D} \bQ \le 0$ for all $\bQ \in
  \mathbb{R}^8$
\item For $\bQ \in \mathbb{R}^8$, the following are equivalent:
  \begin{enumerate}
  \item $\B{D} \bQ = 0$
  \item $(\Phi'(\bQ))^T \B{D} \bQ = 0$
  \end{enumerate}
\item $\Phi''(\bQ)$ is positive-definite
\end{enumerate}
Here the primes indicate gradients with respect to $\bQ$, so $\Phi''$
is the Hessian of $\Phi$ with respect to $\bQ$.  The definition of
Chen, Levermore, and Liu includes an additional clause in item 3
related to an operator we call $\B{\Pi}$ ($\mathcal{Q}$ in their
notation) that maps from $\bQ$ to the conserved quantities of the
relaxation part of the system, $\pr_t \bQ = \B{D} \bQ$ --- namely,
that conditions 3(a) and 3(b) should also be equivalent to
$\Phi'(\bQ)^T = \B{v}^T \B{\Pi}$ for some appropriately-sized vector
$\B{v}$.  Rather than take this as part of the definition of a
strictly convex entropy function, it is more convenient here to take
it as a requirement on $\B{\Pi}$.

Suppose $\Phi(\bQ) = \mathcal{E}(\bQ) = \frac{1}{2} \bQ^T \B{E} \bQ$.
From the preceding sections we already know that conditions 1 and 4
are satisfied, so it only remains to prove conditions 2 and 3.  Since
$\Phi'(\bQ) = \B{E} \bQ$, condition 2 reduces to $\bQ^T \B{E}
\B{D} \bQ \le 0$ for all $\bQ \in \mathbb{R}^8$.  Using equations
\eqref{eq:matrix-D} and \eqref{eq:energy-hessian} for $\B{D}$ and
$\B{E}$, and recalling $\Delta_i = \rho m_i - \rho_f^2$, we find that
\begin{equation} \label{eq:matrix-ed}
  \B{E} \B{D} = \begin{pmatrix}
    0 & 0 & 0 & 0 & 0 & 0 & 0 & 0 \\
    0 & 0 & 0 & 0 & 0 & 0 & 0 & 0 \\
    0 & 0 & 0 & 0 & 0 & 0 & 0 & 0 \\
    0 & 0 & 0 & 0 & 0 & 0 & 0 & 0 \\
    0 & 0 & 0 & 0 & 0 & 0 & 0 & 0 \\
    0 & 0 & 0 & 0 & 0 & 0 & 0 & 0 \\
    0 & 0 & 0 & 0 & 0 & 0 & -\eta/\kappa_1 & 0 \\
    0 & 0 & 0 & 0 & 0 & 0 & 0 & -\eta/\kappa_3 \\
  \end{pmatrix}.
\end{equation}
By inspection, $\B{E} \B{D}$ is a symmetric negative-semidefinite
matrix, so $\Phi'(\bQ)^T \B{D} \bQ \le 0$ for all $\bQ \in
\mathbb{R}^8$ and condition 2 is satisfied.

For condition 3, note that 3(a) implies 3(b) since if $\B{D} \bQ = 0$,
necessarily $\Phi'(\bQ)^T \B{D} \bQ = 0$.  To see that
$\Phi'(\bQ)^T \B{D} \bQ = \bQ^T \B{E} \B{D} \bQ = 0$
implies $\B{D} \bQ = 0$, note that from \eqref{eq:matrix-ed}, we have
$\bQ^T \B{E} \B{D} \bQ = -\frac{\eta}{\kappa_1} q_x^2 -
\frac{\eta}{\kappa_3} q_z^2 = 0$ if and only if $q_x = q_z = 0$.
Since $\B{D}$ only has nonzero entries in the columns corresponding to
$q_x$ and $q_z$, $\Phi'(\bQ)^T \B{D} \bQ = 0$ if and only
if $\B{D} \bQ = 0$.  Therefore condition 3 holds, and $\mathcal{E}$ is
a strictly convex entropy function for the poroelastic system
\eqref{eq:first-order-system}.

\section{Finite volume solution method}
\label{sec:fvm-poro}

\subsection{Wave propagation}

We solve the equations of poroelasticity using a Cartesian grid finite
volume approach.  This section describes the basics of the finite volume
method used here, as well as specifics of how we apply this method
to poroelasticity.  For a comprehensive discussion of this class of
finite volume methods, see LeVeque's book~\cite{rjl:fvm-book}.

The class of finite volume method we use here updates cell averages at
every step by solving a Riemann problem between each pair of adjacent
grid cells.  Thinking of one cell as the ``left'' cell of the problem,
and the other as the ``right'' cell, the Riemann solution process
produces the left-going and right-going fluctuations $\mathcal{A}^-
\Delta \bQ$ and $\mathcal{A}^+ \Delta \bQ$ --- the changes in the cell
variables $\bQ$ caused by the left-going and right-going waves --- along
with a set of waves $\mathcal{W}_i$ with speeds $s_i$ that are used to
implement higher-order correction terms.  With these correction
terms included, the methods used here are second-order accurate.
Where solutions are not smooth, wave limiters can be used on the
higher-order terms to prevent spurious oscillations.  While limiters
can reduce the asymptotic order of accuracy of the solution, they
often decrease the actual value of the error, depending on the norm being
used to measure it and on the grid resolution.  They can also improve
the qualitative behavior of the solution by suppressing dispersive
errors, leading to improved estimates of quantities such as wave
arrival times, and keeping total variation from increasing.

For a homogeneous first-order hyperbolic system such as
\eqref{eq:first-order-hyperbolic}, the left-going and right-going
fluctuations are related to the waves and wave speeds by
\begin{equation} \label{eq:fluctuations}
  \mathcal{A}^+ \Delta \bQ = \sum_{s_i > 0} s_i \mathcal{W}_i, \quad
  \mathcal{A}^- \Delta \bQ = \sum_{s_i < 0} s_i \mathcal{W}_i .
\end{equation}
For a linear problem, such as linear poroelasticity, the waves are
simply eigenvectors of the flux Jacobian matrix (for instance, the
matrix $\bA$ of \eqref{eq:matrix-A} for waves propagating in the $x$
direction in a material with its principal axes aligned with the
coodinate axes) associated with the material through which the wave
propagates --- that is, they have the form $\mathcal{W}_i = \beta_i
\B{r}_i$, where $\B{r}_i$ is the eigenvector and $\beta_i$ is a
scalar that gives the strength of the wave.  Each wave speed $s_i$ is
the corresponding eigenvalue of the flux Jacobian.

A quantity of critical importance in these solution methods is the
{\em CFL number} $\nu$.  Informally, the CFL number is the ratio of
the distance a wave travels in one timestep to the width of a grid
cell; more formally, for a Cartesian grid the global CFL number is
\begin{equation} \label{eq:cfl-def}
  \nu = \max_{\text{all cells, all waves}} \max \left( \frac{|s_x|
    \Delta t}{\Delta x}, \frac{|s_z| \Delta t}{\Delta z} \right).
\end{equation}
Here $s_x$ and $s_z$ are the speeds of waves generated from the
Riemann problems in the $x$ and $z$ directions, $\Delta x$ and $\Delta
z$ are the grid spacings, and $\Delta t$ is the time
step size.  The methods used here are stable for $\nu \le 1$; since $\nu$
comes from a maximum over all waves, this means that our stability is
limited by the fast P wave.

Because the poroelasticity equations we use here are a linear system,
solution of the Riemann problem is straightforward.  There is one
complication, however.  We wish to consider domains composed of
multiple materials --- in
fact, our code has the capability for each grid cell to be made of a
different material --- so the coefficient matrices are only piecewise
constant.  We always choose the grid boundaries to coincide with the
material boundaries, but we must still solve Riemann problems between
domains with different coefficient matrices $\bA$ and $\bB$.  Suppose we
are solving a Riemann problem in the $x$-direction, with $\bA = \bA_l$ in
the left cell and $\bA = \bA_r$ in the right cell.  We know that in the
left cell, the Riemann solution consists of waves with strength
$\beta_{li}$ in the directions of the eigenvectors $\B{r}_{li}$ of $\bA_l$,
corresponding to the negative eigenvalues of $\bA_l$; similarly, in the
right cell we will have waves with strength $\beta_{ri}$ in the
directions of eigenvectors $\B{r}_{ri}$ or $\bA_r$, corresponding to the
positive eigenvalues of $\bA_r$.  There will also be a stationary
discontinuity at the cell interface, which will lie in the null space
of $\bA_l$ and $\bA_r$.  (The fact that $\bA$ has the same null space
for any poroelastic material greatly simplfies matters here, and the corresponding
eigenvectors $\B{r}_4$ and $\B{r}_5$ will not carry a subscript identifying
them with the left or right material.)  The total jump in $\bQ$ across
all the waves and the stationary discontinuity must add up to the
difference in $\bQ$ between the left and right states, $\Delta \bQ = \bQ_r -
\bQ_l$, so we require
\begin{equation}
  \Delta \bQ =
  \sum_{i=1}^3 \beta_{li}\B{r}_{li} + \sum_{i=4}^5 \beta_i\B{r}_i +
  \sum_{i=6}^8 \beta_{ri}\B{r}_{ri} =: \tilde{\B{R}} \bbeta .
\end{equation}
We can thus compute the wave strengths as $\bbeta = \tilde{\B{R}}^{-1}
\Delta \bQ$.  In practice, since the strength of the stationary
discontinuity is never used directly, we never compute $\beta_4$ and
$\beta_5$.  The same analysis holds for a Riemann problem in the
$z$-direction.  This approach corresponds to an open-pore condition
between the two poroelastic media, as described by Deresiewicz and
Skalak~\cite{deresiewicz-skalak:poro-uniqueness-1963} and validated by
Gurevich and Schoenberg~\cite{gurevich-schoenberg:interface-1999}.
Other interface conditions are possible in poroelasticity, such as the
closed or partially open pore conditions of Deresiewicz and Skalak, or
the loose contact condition of
Sharma~\cite{sharma:dissimilar-boundary-2008}; we do not model these
in this work, but they would be straightforward to incorporate into
the Riemann solution process.

Because the eigenstructure of poroelasticity is somewhat complex, we
do not compute the eigensystems of $\bA$ and $\bB$ analytically.
Instead, we use LAPACK~\cite{lapack-users-guide} to compute the
eigenvalues and eigenvectors of $\bA$ and $\bB$ for every poroelastic
material present in the model, and the $\tilde{\B{R}}^{-1}$ matrices
for each Riemann solve direction and every pair of left and right
materials that could occur.  For efficiency, we pre-compute these
quantities for all materials used (or all possible pairs of
materials in the case of $\tilde{\B{R}}^{-1}$) before starting the
solution proper, and look them up using a material number stored
with each cell during the Riemann solves.

Boundary conditions were implemented using the usual ghost cell
approach~\cite{rjl:fvm-book}.  We set ghost cell values using either
zero-order extrapolation, for boundaries where waves should flow
outward and not return, or by setting the ghost cell values equal to
the exact solution at the centers of those cells, when we verified our
code against known analytic solutions.

\subsection{Operator splitting}

We include the dissipative part $\bQ_t = \B{D}\bQ$ of the poroelasticity
equations using operator splitting.  Since the $\B{D}$ matrix is constant,
we can use the exact solution operator $\exp(\B{D}\Delta t)$ to advance the
solution by a time increment $\Delta t$; not only is this the most
accurate solution available for this part of the system, it is also
unconditionally stable and allows the time step to be chosen based
solely on stability for the wave propagation part of the system.

The software framework we use offers either Godunov or Strang
splitting as a run-time option.  Godunov splitting is formally
first-order accurate in time and uses a single full-length step of the
source term operator per time step, while Strang splitting is
second-order and uses two half-steps of the source term.  For many
practical problems Godunov splitting is a good choice because it
displays similar similar error to Strang --- the coefficient of the
first-order error term is often small --- while being less
computationally intensive.  However, we primarily use Strang splitting
here because it displays substantially greater accuracy for the
particular poroelasticity problems we solve, and because our source
term is computationally cheap compared to the wave propagation part of
the system.  For comparison, we also show results for Godunov
splitting.

\subsection{Stiff regime and subcharacteristic condition}
\label{sec:stiff-subcharacteristic}

For some cases we consider, the time step is much larger
than the characteristic time scales associated with the solution of
$\pr_t \bQ = \bD\bQ$.  These cases fall outside the regime where
asymptotic leading-order error
estimates
are relevant, and for them the source
term is stiff, a known source of difficulty in operator splitting
approaches for hyperbolic
equations~\cite{colella-majda-roytburd:reacting-shock,
  rjl-yee:stiff-split}.  Based on a conjecture of
Pember~\cite{pember:spurious-solns}, we expect to avoid spurious
solutions for this stiff relaxation system if the poroelasticity
equations \eqref{eq:first-order-system} satisfy a {\em
  subcharacteristic condition}, where the wave speeds for the reduced
equations obtained by restricting the full system to the equilibrium
manifold of the dissipation term (i.e. zero fluid velocity relative to
the solid matrix) interlace with the wave speeds for the full system.
The appropriate generalization of Pember's conjecture to systems of
more than two equations is not obvious, but based on the principle
that information should propagate more slowly (certainly no more
quickly!)\ in the reduced system than in the full system, we expect to
avoid spurious solutions if for all possible wave propagation
directions the speeds $\lambda_1$ and $\lambda_2$ of the reduced
system are strictly less than the speed of the fastest wave of the
full system,
\begin{equation} \label{eq:subcharacteristic}
  \lambda_1 < c_{pf}, \quad \lambda_2 < c_{pf}.
\end{equation}
Here $c_{pf}$ is the speed of a fast P wave.  We ignore the negative
eigenvalues for both the full and reduced systems, because they are
simply the negatives of the wave speeds and will automatically satisfy
a similar inequality.  We also ignore the zero eigenvalues of the full
system, since they correspond to eigencomponents of the solution that
are left unchanged in the wave propagation part of the solution
process, and are evolved according to the exact solution operator
$\exp(\bD\Delta t)$ in the dissipation part.

To construct the appropriate reduced system, we follow the derivation
of Chen, Levermore, and
Liu~\cite{chen-levermore-liu:stiff-relaxation}.  First we examine the
dissipation part of the system in isolation,
\begin{equation} \label{eq:reduced-sys-dissipation}
  \pr_t \bQ = \bD\bQ.
\end{equation}
System \eqref{eq:reduced-sys-dissipation} has six conserved quantities
$\B{u} = \B{\Pi} \bQ$, which are related to the state variables $\bQ$ by the
matrix
\begin{equation} \label{eq:matrix-pi}
  \B{\Pi} := \begin{bmatrix}
    1 & 0 & 0 & 0 & 0 & 0 & 0 & 0 \\
    0 & 1 & 0 & 0 & 0 & 0 & 0 & 0 \\
    0 & 0 & 1 & 0 & 0 & 0 & 0 & 0 \\
    0 & 0 & 0 & 1 & 0 & 0 & \rho_f/\rho & 0 \\
    0 & 0 & 0 & 0 & 1 & 0 & 0 & \rho_f/\rho \\
    0 & 0 & 0 & 0 & 0 & 1 & 0 & 0
  \end{bmatrix}.
\end{equation}
The fact that $\B{u}$ is a vector of conserved quantities of
\eqref{eq:reduced-sys-dissipation} follows immediately from the fact
that $\B{\Pi} \bD = 0$, so that $\pr_t \B{u} = \B{\Pi} \pr_t \bQ =
\B{\Pi} \bD \bQ = 0$.  Given the conserved quantities $\B{u}$, the
unique equilibrium $\bQ_{\text{eq}}$ of
\eqref{eq:reduced-sys-dissipation} that satisfies $\bD \bQ_{\text{eq}}
= 0$ and $\B{\Pi}\bQ_{\text{eq}} = \B{u}$ is $\bQ_{\text{eq}} =
\bG\B{u}$, where
\begin{equation}
  \bG := \begin{bmatrix}
    \B{I}_{6 \times 6} \\ \B{0}_{2 \times 6}
  \end{bmatrix}.
\end{equation}
Notice that $\B{\Pi}$ and $\bG$ satisfy the relation $\B{\Pi} \bG
= \B{I}_{6\times 6}$.  The reduced system is found by multiplying the full
poroelastic system from the left by $\B{\Pi}$ to eliminate the
dissipation term, and requiring the state vector $\bQ$ to lie on the
equilibrium manifold, $\bQ = \bQ_{\text{eq}} = \bG\B{u}$, resulting in
\begin{equation} \label{eq:reduced-sys}
  \pr_t \B{u}_t + \B{\Pi} \bA\bG \pr_x \B{u} + \B{\Pi} \bB\bG \pr_z \B{u} = 0.
\end{equation}

Now that we have the matrix $\B{\Pi}$, we can show that the additional
condition of Chen, Levermore, and Liu mentioned in Section
\ref{sec:entropy-function} also holds --- that the statements $\bD \bQ
= 0$ and $\Phi'(\bQ)^T \B{D} \bQ = \bQ^T \B{E} \B{D} \bQ =
0$ are equivalent to $\Phi'(\bQ)^T = \bQ^T \B{E} = \B{v}^T
\B{\Pi}$ for some $\B{v} \in \mathbb{R}^6$.  First, if $\bQ^T \B{E} =
\B{v}^T \B{\Pi}$, we immediately have $\bQ^T \B{E} \bD \bQ = \B{v}^T
\B{\Pi} \bD \bQ = 0$ since $\B{\Pi} \bD = 0$.  We also noted in
section \ref{sec:entropy-function} that if $\bD \bQ = 0$, then $q_x = q_z =
0$.  Thus the fourth, fifth, seventh, and eighth components of $\bQ^T
\B{E}$ are $(\bQ^T \B{E})_4 = \rho v_x$, $(\bQ^T \B{E})_5 = \rho v_z$,
$(\bQ^T \B{E})_7 = \rho_f v_x$, and $(\bQ^T \B{E})_8 = \rho_f v_z$.
We therefore have $\bQ^T \B{E} = \B{v}^T \B{\Pi}$ with $\B{v}$ given by
\begin{equation}
 \B{v} = \begin{pmatrix}
   (\bQ^T \B{E})_1 \\
   (\bQ^T \B{E})_2 \\
   (\bQ^T \B{E})_3 \\
   \rho v_x \\
   \rho v_z \\
   (\bQ^T \B{E})_6
 \end{pmatrix}.
\end{equation}

Given that a strictly convex entropy function $\Phi(\bQ) =
\mathcal{E}(\bQ)$ exists and the above condition is satisfied, Theorem
2.1 and the subsequent remark of Chen, Levermore, and
Liu~\cite{chen-levermore-liu:stiff-relaxation} imply that the reduced
system \eqref{eq:reduced-sys} is hyperbolic, and satisfies a nonstrict
subcharacteristic condition, which we can render here in the context
of poroelasticity as
\begin{equation} \label{eq:cll-nonstrict-subcharacteristic}
    c_{ps} \le \lambda_2 \le c_{pf}, \quad 0 \le \lambda_1 \le c_{s}.
\end{equation}
While this does not imply the strict subcharacteristic condition
\eqref{eq:subcharacteristic}, it is nearly as useful.  This
characterization of the eigenvalues comes from the fact that
eigenvalues of the symmetric generalized eigenproblem
\eqref{eq:eigenproblem-energy} satisfy the Rayleigh quotient minimax
principle
\begin{equation} \label{eq:eigenvalue-minimax}
    \lambda_{fk}^\downarrow = \max_{\mathcal{S}_k} \min_{\B{v} \in \mathcal{S}_k}
    \frac{\B{v}^T \B{E} \breve{\bA} \B{v}}{\B{v}^T \B{E} \B{v}}, \quad
    \lambda_{fk}^\uparrow = \min_{\mathcal{S}_k} \max_{\B{v} \in \mathcal{S}_k}
    \frac{\B{v}^T \B{E} \breve{\bA} \B{v}}{\B{v}^T \B{E} \B{v}},
\end{equation}
where $\mathcal{S}_k$ is any $k$-dimensional subspace of
$\mathbb{R}^8$, $\lambda_{fk}^\downarrow$ is the $k$'th eigenvalue of
the full system counting down from the largest, and
$\lambda_{fk}^\uparrow$ is the $k$'th eigenvalue counting up from the
smallest.  Chen, Levermore, and Liu prove that the eigenvalues of the
reduced system satisfy a similar minimax principle over a restricted
set of subspaces,
\begin{equation} \label{eq:eigenvalue-minimax-restricted}
  \begin{aligned}
    \lambda_{ek}^\downarrow = \max_{\mathcal{S}_k \subseteq \mathcal{R}(\bG)} \min_{\B{v} \in \mathcal{S}_k}
    \frac{\B{v}^T \B{E} \breve{\bA} \B{v}}{\B{v}^T \B{E} \B{v}}, \quad
    \lambda_{ek}^\uparrow = \min_{\mathcal{S}_k \subseteq \mathcal{R}(\bG)} \max_{\B{v} \in \mathcal{S}_k}
    \frac{\B{v}^T \B{E} \breve{\bA} \B{v}}{\B{v}^T \B{E} \B{v}},
  \end{aligned}
\end{equation}
where the subscript $e$ indicates eigenvalues of the reduced
(``equilibrium'') system and $\mathcal{R}(\bG)$ is the range space of
$\bG$.

Equality can be realized in the nonstrict subcharacteristic condition
\eqref{eq:cll-nonstrict-subcharacteristic} --- for example, if the
fluid density $\rho_f$ for the orthotropic sandstone material whose
parameters are given in the first column of Table \ref{tab:matprops}
is reduced from $1040$ kg/m$^3$ to $208.9$ kg/m$^3$, a fast P wave
traveling along the material principal $1$-axis shows no fluid
relative motion.  This means its eigenvector lies in
$\mathcal{R}(\bG)$, so by \eqref{eq:eigenvalue-minimax-restricted} the
reduced system shares the same wave speed.  Fortunately, this turns
out to be innocuous from the standpoint of the true solution --- if
the eigenvector $\B{v}$ associated with some wave is in
$\mathcal{R}(\bG)$, then it is in $\mathcal{N}(\bD)$;
since it is an eigenvector of both parts of the system, the
corresponding eigencomponent of the solution can be decoupled from the
rest of the system, and its solution is independent of them.
Furthermore, the PDE describing the evolution of this component of the
solution is a purely hyperbolic one, with no source term --- there is
no momentum transfer due to viscous drag between the solid and fluid
for this wave mode because there is no relative motion between them.
A similar result carries over to the numerical solution obtained by
operator splitting: if $\B{v}$ is in $\mathcal{N}(\bD)$, then
$\exp(\bD \Delta t) \B{v} = \B{v}$, and the corresponding component of
$\bQ$ passes through the solution operator for the dissipation term
unchanged.

As an aside, the reduced system \eqref{eq:reduced-sys} has a familiar
form.  If we multiply out the coefficient matrices, we get
\begin{equation} \label{eq:reduced-sys-matrices}
  \B{\Pi} \bA \bG = \begin{pmatrix}
    0 & 0 & 0 & -c^u_{11} & 0 & 0 \\
    0 & 0 & 0 & -c^u_{13} & 0 & 0 \\
    0 & 0 & 0 & 0 & -c^u_{55} & 0 \\
    -\frac{1}{\rho} & 0 & 0 & 0 & 0 & 0 \\
    0 & 0 & -\frac{1}{\rho} & 0 & 0 & 0 \\
    0 & 0 & 0 & M \alpha_{1} & 0 & 0
  \end{pmatrix}, \quad
  \B{\Pi} \bB \bG = \begin{pmatrix}
    0 & 0 & 0 & 0 & -c^u_{13} & 0 \\
    0 & 0 & 0 & 0 & -c^u_{33} & 0 \\
    0 & 0 & 0 & -c^u_{55} & 0 & 0 \\
    0 & 0 & -\frac{1}{\rho} & 0 & 0 & 0 \\
    0 & -\frac{1}{\rho} & 0 & 0 & 0 & 0 \\
    0 & 0 & 0 & 0 & M \alpha_{3} & 0
  \end{pmatrix}.
\end{equation}
The upper-left $5\times 5$ portions of these matrices are just the
coefficient matrices for orthotropic plane-strain elasticity, with the
fluid pressure coming along as an additional variable
determined entirely by the elastic field variables.  Because of this
we will identify the faster wave of the reduced system as the
``reduced P wave'' and the slower one as the ``reduced S wave.''
Note that for this work, the reduced system is only of theoretical
importance --- for the actual numerical code, we discretize the full
system \eqref{eq:first-order-system}.

\subsection{Numerical software}

We implemented the numerical solution techniques described here using
the \clawpack{} finite volume method package, version
4.6~\cite{claw.org}.  \clawpack{} implements the parts of a
high-resolution finite volume code that are common across all
problems, leaving the user to write only problem-specific code such as
Riemann solvers.  Operator splitting is supported for source terms,
such as the dissipative term here, by means of a user-supplied
subroutine that advances the system by a specified time step under the
action of the source term.  Both Godunov and Strang splitting are
available.  Block-structured Berger-Colella-Oliger adaptive mesh
refinement (AMR) is available from the \amrclaw{}
package~\cite{mjb-rjl:amrclaw}; \amrclaw{} can also run in parallel on
shared-memory systems using OpenMP.  Besides Cartesian grids,
\clawpack{} and \amrclaw{} also support logically rectangular mapped
grids --- while we do not use mapped grids here, we plan to present
results with them in a subsequent publication.

\section{Results}
\label{sec:results}

We present results here for four classes of problems.  First, we
demonstrate convergence of our numerical solution to known analytic
plane wave solutions for an orthotropic medium for the wave
propagation part of the system alone.  We then include the viscous
dissipation term, and examine the effect of operator splitting on
accuracy, again comparing against known analytic plane wave solutions.
Next, we show results for simple point sources in uniform othrotropic
media, solving test problems previously addressed by de la Puente et
al.~\cite{delapuente-dumbser-kaser-igel:poro-dg} and
Carcione~\cite{carcione:poro-aniso-1996} in order to further verify
our code.  Finally, we solve a larger-scale problem with a domain
composed of two isotropic materials, involving wave reflection,
refraction, and interconversion at the material boundary, as well as
demonstrating the use of adaptive mesh refinement to reduce the time
required for solution.

\subsection{Analytic plane wave solution}

Before we test our code's convergence against analytic plane wave
solutions, however, we will first outline the procedure used to obtain these
analytic solutions.  We start by assuming the velocity and stress
fields have a plane wave form,
\begin{alignat}{4}
  \B{V} &:= \begin{pmatrix} v_x & v_z & q_x & q_z \end{pmatrix}^T &&=
  \B{V}_0 \exp(i(k_x x + k_z z - \omega t))\\
  \B{T} &:= \begin{pmatrix} \tau_{xx} & \tau_{zz} & \tau_{xz} &
    -p \end{pmatrix}^T &&= \B{T}_0 \exp(i(k_x x + k_z z - \omega t)).
\end{alignat}
Here $\B{V}_0$ and $\B{T}_0$ are constant vectors, and $\omega$ is the
prescribed angular frequency of the wave.  The wavenumbers $k_x$ and
$k_z$ are yet to be determined, but we also prescribe $k_x = k l_x$
and $k_z = k l_z$, where $l_x$ and $l_z$ are the (real-valued)
direction cosines of the wavevector, with $l_x^2 + l_z^2 = 1$.

With these assumptions on the solution, stress-strain equations
\eqref{stress-strain-1} through \eqref{stress-strain-4} imply
\begin{equation} \label{eq:plane-constitutive}
  -\omega \B{T}_0 = k\B{F} \B{V}_0,
\end{equation}
where the matrix $\B{F}$ is
\begin{equation}
  \B{F} = \begin{pmatrix}
    l_x c_{11}^u & l_z c_{13}^u & \alpha_1 M l_x & \alpha_1 M l_z\\
    l_x c_{13}^u & l_z c_{33}^u & \alpha_3 M l_x & \alpha_3 M l_z \\
    l_z c_{55}^u & l_x c_{55}^u & 0 & 0\\
    \alpha_1 M l_x  & \alpha_3 M l_z &  M l_x & M l_z
  \end{pmatrix}.
\end{equation}
Equations of motion \eqref{motion-x} through \eqref{darcy-z} also
imply
\begin{equation} \label{eq:plane-eom}
  k \B{L} \B{T}_0 = -\omega \bGamma \B{V}_0,
\end{equation}
where the matrices $\B{L}$ and $\bGamma$ are
\begin{equation}
  \B{L} = \begin{pmatrix}
    l_x & 0 & l_z & 0 \\
    0 & l_z & l_x & 0\\
    0 & 0 & 0 & l_x \\
    0 & 0 & 0 & l_z
  \end{pmatrix}, \quad
  \bGamma = \begin{pmatrix}
    \rho & 0 & \rho_f & 0\\
    0 & \rho & 0 &\rho_f\\
    \rho_f & 0 & iY_1(-\omega)/\omega & 0\\
    0 & \rho_f & 0 & iY_3(-\omega)/\omega 
  \end{pmatrix},
\end{equation}
and $Y_j(\omega) := i \omega m_j + \eta/\kappa_j$ for $j = 1,3$.

Combining equations \eqref{eq:plane-constitutive} and
\eqref{eq:plane-eom}, we can obtain an eigenproblem for $\B{V}_0$ and
$\left(\frac{\omega}{k}\right)^2$:
\begin{equation} \label{eq:plane-eigenproblem}
  \bGamma^{-1} \B{L} \B{F} \B{V}_0 = \left(\frac{\omega}{k}\right)^2 \B{V}_0.
\end{equation}
To obtain a plane wave solution, we solve this eigenproblem, choose
the eigenvalue $\left(\frac{\omega}{k}\right)^2$ and eigenvector
$\B{V}_0$ corresponding to the wave family of interest, then back out
the appropriate $\B{T}_0$ from \eqref{eq:plane-constitutive}.  Note
that the eigenvalues and eigenvectors will be complex-valued if
dissipation is present.

\subsection{Plane wave convergence results --- inviscid}

We conducted our plane wave convergence tests on a uniform domain
composed of orthotropic layered sandstone, whose properties are given in Table
\ref{tab:matprops}.  The wave-propagation part of our code was first
tested alone, without viscous dissipation included.  These tests were
conducted using plane waves with a fixed angular frequency of $10^4$
rad/s, in a square domain $8\,\text{m}$ on a side.  This distance is about
two wavelengths of the fast P wave at this frequency.  The total
simulation time was $2\pi \times 10^{-4}\,\text{s}$ --- one period of the
wave.  Because there is no intrinsic time scale associated with Biot's
equations when viscous dissipation is omitted, the frequency of the
plane wave is not directly relevant to the accuracy of these tests;
instead, only the number of grid cells per wavelength and the ratio of
the simulation time to the wave period are relevant.  Boundary
conditions were set by filling the ghost cells with the value of the
true plane wave solution at the cell centers, and the initial
condition on the grid was set the same way.  This ensured that the
accuracy of the numerical solution was governed only by the
correctness of the wave-propagation algorithm used within the problem
domain, not by the implementation of the boundary conditions.

\begin{table}[ht]
\caption{Properties of the poroelastic media used in test cases, taken
  from de la Puente et al.~\cite{delapuente-dumbser-kaser-igel:poro-dg}.  Wave speeds are
  correct in the high-frequency limit.  $c_{pf}$ is the fast P wave
  speed, $c_s$ is the S wave speed, $c_{ps}$ is the slow P wave speed,
  and $\tau_d$ is the time constant for dissipation.  Subscript numbers
  indicated principal directions.}
\label{tab:matprops}
\begin{center}
\begin{tabular}{rlp{0.8in}p{0.8in}p{0.8in}p{0.8in}p{0.8in}}
  & & Sandstone (orthotropic) & Glass/epoxy (orthotropic) & Sandstone
  (isotropic) & Shale (isotropic) \\
  \hline
  \multicolumn{6}{l}{Base properties} \\
  \hline
  $K_s$ & (GPa) & 80 & 40 & 40 & 7.6 \\
  $\rho_s$ & (kg/m$^3$) & 2500 & 1815 & 2500 & 2210 \\
  $c_{11}$ & (GPa) & 71.8 & 39.4 & 36 & 11.9 \\
  $c_{12}$ & (GPa) & 3.2 & 1.2 & 12 & 3.96 \\
  $c_{13}$ & (GPa) & 1.2 & 1.2 & 12 & 3.96 \\
  $c_{33}$ & (GPa) & 53.4 & 13.1 & 36 & 11.9 \\
  $c_{55}$ & (GPa) & 26.1 & 3.0 & 12 & 3.96 \\
  $\phi$ & & 0.2 & 0.2 & 0.2 & 0.16 \\
  $\kappa_1$ & ($10^{-15}$ m$^2$) & 600 & 600 & 600 & 100 \\
  $\kappa_3$ & ($10^{-15}$ m$^2$) & 100 & 100 & 600 & 100 \\
  $T_1$ & & 2 & 2 & 2 & 2 \\
  $T_3$ & & 3.6 & 3.6 & 2 & 2 \\
  $K_f$ & (GPa) & 2.5 & 2.5 & 2.5 & 2.5 \\
  $\rho_f$ & (kg/m$^3$) & 1040 & 1040 & 1040 & 1040 \\
  $\eta$ & ($10^{-3}$ kg/m$\cdot$s) & 1 & 1 & 0 & 0 \\ \hline
  \multicolumn{6}{l}{Derived quantites} \\
  \hline
  $c_{pf1}$ & (m/s) & 6000 & 5240 & 4250 & 2480 \\
  $c_{pf3}$ & (m/s) & 5260 & 3580 & 4250 & 2480 \\
  $c_{s1}$ & (m/s) & 3480 & 1370 & 2390 & 1430 \\
  $c_{s3}$ & (m/s) & 3520 & 1390 & 2390 & 1430 \\
  $c_{ps1}$ & (m/s) & 1030 & 975 & 1020 & 1130 \\
  $c_{ps3}$ & (m/s) & 746 & 604 & 1020 & 1130 \\
  $\tau_{d1}$ & ($\mu$s) & 5.95 & 5.85 & --- & ---\\
  $\tau_{d3}$ & ($\mu$s) & 1.82 & 1.81 & --- & ---
\end{tabular}
\end{center}
\end{table}

\begin{figure}
  \begin{center}
    \includegraphics[width=0.5\textwidth]{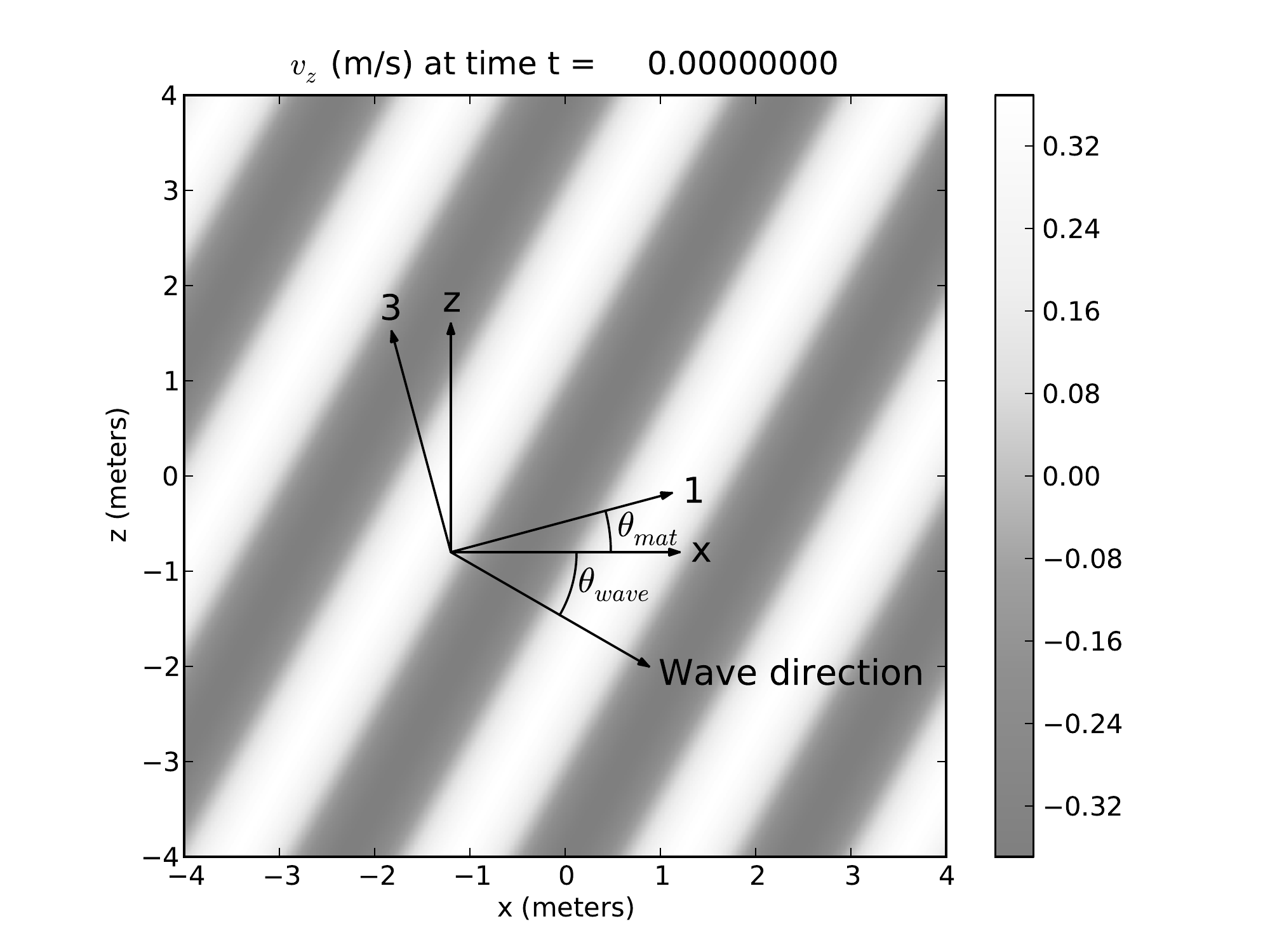}
    \caption{Sample plane wave initial condition, showing the relation
      between the global $x-z$ axes, the material $1-3$ axes, and the
      wave propagation direction.  For this plot $\theta_{\text{mat}}
      = 15^\circ$ and $\theta_{\text{wave}} = -30^\circ$.  The shading
      shows the $z$-direction solid matrix velocity.  This plot depicts
      a fast P wave, with an angular frequency of $10^4$ rad/s, in the
      same orthotropic sandstone medium used for the inviscid
      convergence tests.  The grid dimensions are $800 \times 800$
      cells.  \label{fig:plane-wave-diagram}}
  \end{center}
\end{figure}

Since we are working with an orthotropic material, rather than an
isotropic one, the speed and associated eigenvector for a plane wave
depend on its propagation direction --- the solutions for plane waves
propagating in different directions are not simply rotated versions of
each other, and in order to be confident in the correctness of our
code it was necessary to test it with plane waves propagating at a
variety of angles $\theta_{\text{wave}}$ relative to the global $x$
axis.  Since we anticipate solving problems where the principal
material axes do not coincide with the global coordinate axes, we also
tested our code with a variety of angles $\theta_{\text{mat}}$ between
the material 1 axis and the $x$ axis.  Figure
\ref{fig:plane-wave-diagram} shows a sample plane wave solution, with
the relevant axes and angles identified.  For the studies presented
here, we used $\theta_{\text{wave}}$ values from $0^\circ$ to
$345^\circ$ counterclockwise, and $\theta_{\text{mat}}$ values from
$0^\circ$ to $165^\circ$ counterclockwise, both in steps of
$15^\circ$.  (When expressed in global $x$-$z$ coordinates, the system
matrices $\bA$ and $\bB$ simply flip sign after a rotation of
$\theta_{\text{mat}} = 180^\circ$, while $\bD$ is unchanged, so it was
not necessary to use $\theta_{\text{mat}}$ values of $180^\circ$ or
over.)  This gives 24 different $\theta_{\text{wave}}$ values and 12
different $\theta_{\text{mat}}$ values, for a total of 288
combinations of these angles.  For each $(\theta_{\text{wave}},
\theta_{\text{mat}})$ pair, we examined plane waves in each of the
three families, on $100 \times 100$, $200 \times 200$, $400 \times
400$, and $800 \times 800$ cell grids, for a total of 3456 different
test cases.  For all cases, we chose the time step so that the CFL
number was 0.9.  Because the solution is smooth, we used no wave
limiting for this convergence study.

For each test case, we measured the error by taking the energy norm
on each cell of the difference between the numerically obtained cell
value and the true solution at the cell center.  We then applied
the grid 1-, 2-, and max-norms to the energy norm error field to
obtain an aggregate error norm.  (The grid 1-norm used here is just
the linear algebraic 1-norm divided by the number of grid cells;
similarly, the grid 2-norm is the linear algebraic 2-norm divided by
the square root of the number of grid cells.)  For each combination of
$\theta_{\text{wave}}$, $\theta_{\text{mat}}$, and wave family, we
performed a linear least-squares fit of $\log(\text{error})$ versus
$\log(m)$, where $m$ is the number of grid cells along each axis.  The
slope of this fit was considered to be the convergence rate of the
code for this set of cases.

Table \ref{tab:conv-summary-inviscid} summarizes the results of this
convergence study.  The convergence rates listed are the maximum,
minimum, and mean over all combinations of $\theta_{\text{wave}}$ and
$\theta_{\text{mat}}$; the worst $R^2$ value for the least-squares fit
of $\log(\text{error})$ versus $\log(m)$ is also reported.  The last
two columns of Table \ref{tab:conv-summary-inviscid} give the lowest
and highest error for each wave in each norm on the finest grid,
observed over all combinations of $\theta_{\text{wave}}$ and
$\theta_{\text{mat}}$, normalized by the energy norm of the plane wave
eigenvector $\bQ_0$ formed from the vectors $\B{V}_0$ and $\B{T}_0$ of
the preceding section in order to allow a fair comparison
between different cases.  We see second-order convergence in all three
norms for the fast P and S waves.  Results are impaired for
the slow P-wave because its slow propagation speed causes it to be
underresolved on the coarser grids at the frequency used, but we still
observe second-order convergence in the 1-norm and 2-norm.
While the error
varies somewhat depending on the wave propagation and principal material
directions, it does so by no more than a factor of 3-4, indicating
that there are no severe grid alignment effects.

\begin{table}
\caption{Summary of convergence results for inviscid test cases}
\label{tab:conv-summary-inviscid}
\begin{center}
\begin{tabular}{cccccccc}
  & & \multicolumn{3}{c}{Convergence rate} & &
  \multicolumn{2}{c}{Error on $800 \times 800$ grid}\\
  \cline{3-5} \cline{7-8}
  Wave family & Error norm & Best & Worst & Mean & Worst $R^2$ value &
  Best & Worst \\
  \hline
  \multirow{3}{*}{Fast P} & 1-norm & 2.03 & 2.01 & 2.01 & 0.99994
  & $2.53\times 10^{-5}$ & $6.90\times 10^{-5}$ \\
  & 2-norm & 2.02 & 2.01 & 2.01 & 0.99995 & $3.04\times 10^{-5}$ &
  $8.23\times 10^{-5}$ \\
  & Max-norm & 2.02 & 1.96 & 2.00 & 0.99943 & $5.61\times 10^{-5}$ &
  $1.76\times 10^{-4}$ \\
  \hline
  \multirow{3}{*}{S} & 1-norm & 2.01 & 2.00 & 2.01 & 1.00000 &
  $1.31\times 10^{-4}$ & $3.23\times 10^{-4}$ \\
  & 2-norm & 2.01 & 2.00 & 2.00 & 1.00000 & $1.49\times 10^{-4}$ &
  $3.74\times 10^{-4}$ \\
  & Max-norm & 2.00 & 1.96 & 1.99 & 0.99977 & $2.80\times 10^{-4}$ &
  $7.98\times 10^{-4}$ \\
  \hline
  \multirow{3}{*}{Slow P} & 1-norm & 1.99 & 1.92 & 1.97 & 0.99922 &
  $3.05\times 10^{-3}$ & $1.13\times 10^{-2}$ \\
  & 2-norm & 1.99 & 1.94 & 1.97 & 0.99953 & $3.44\times 10^{-3}$ &
  $1.27\times 10^{-2}$\\
  & Max-norm & 1.93 & 1.67 & 1.80 & 0.99146 & $8.81\times 10^{-3}$ &
  $3.16\times 10^{-2}$
\end{tabular}
\end{center}
\end{table}

\begin{figure}
  \begin{center}
  \includegraphics[width=0.8\textwidth]{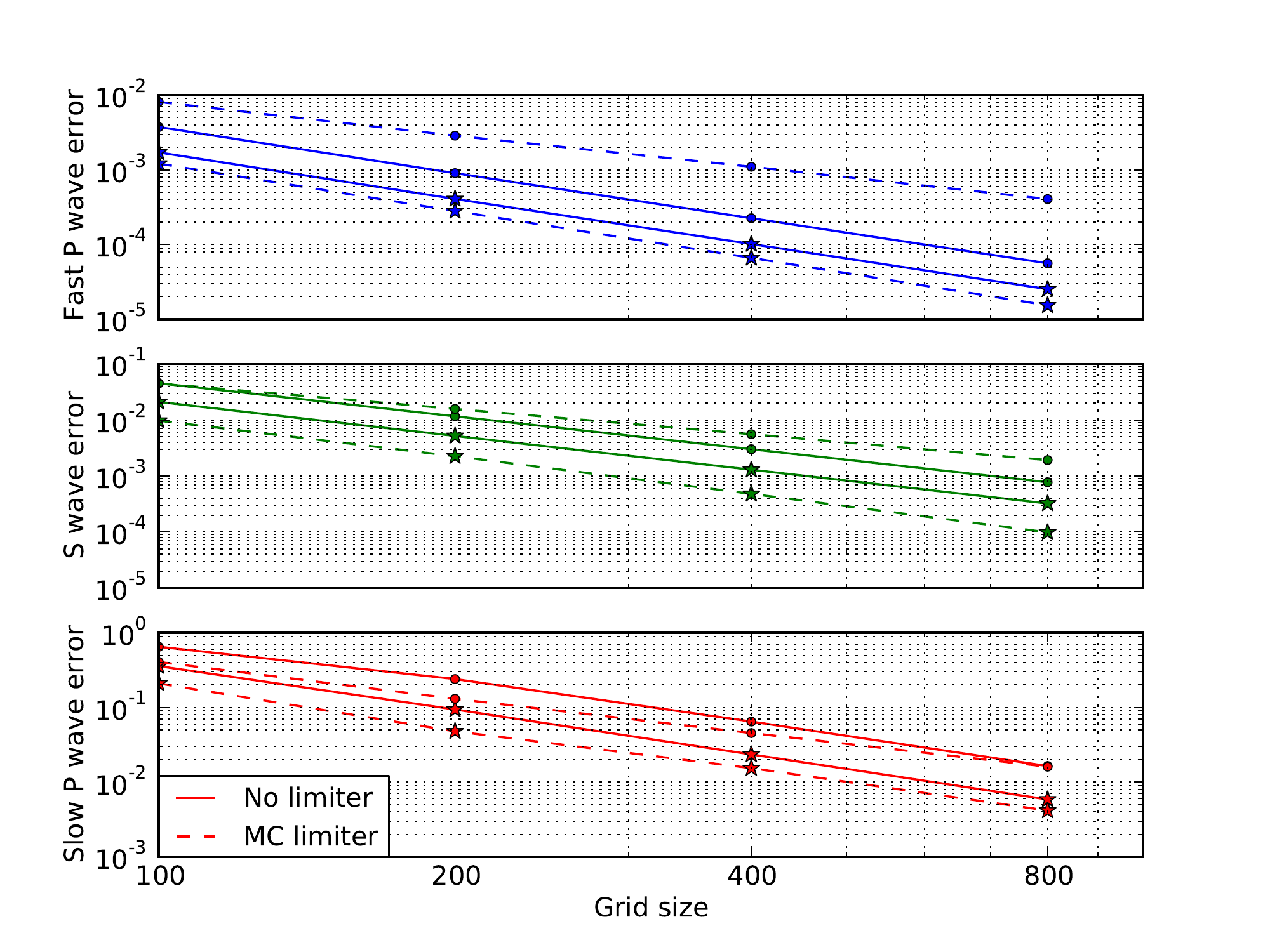}
  \caption{Comparison of error both with and without a limiter for the
    inviscid plane wave test cases with $\theta_{\text{mat}} =
    \theta_{\text{wave}} = 0$.  The circles indicate normalized energy
    max-norm error, while the stars indicate normalized energy 1-norm
    error.  The limiter is always beneficial in the 1-norm for these
    cases, and is also helpful in the max-norm for poorly-resolved
    waves.
    \label{fig:limiter-nolimiter}}
  \end{center}
\end{figure}

We do not include convergence results with wave limiting here, but
informal exploration suggests that for the cases above, limiting
reduces the order of accuracy to about 1.9 in the 1-norm, 1.8 in the
2-norm, and 1.6 in the max-norm.  This is because wave limiting tends
to clip extrema in order to avoid introducing spurious oscillations.
Despite the reduced order of accuracy, however, using a limiter can
improve actual error in many cases, often in the 1-norm but even in
the max-norm if a wave is poorly resolved or heavily affected by
dispersive errors.  Figure \ref{fig:limiter-nolimiter} shows an
example of the effect of limiting, using the Monotonized Centered (MC)
limiter, on the normalized energy max-norm and 1-norm errors in each
of the three wave families for the inviscid test cases with
$\theta_{\text{mat}} = \theta_{\text{wave}} = 0$.  For each wave, the
max-norm error decreases more slowly with increasing grid size when
the limiter is present.  The fast P wave is well resolved even on the
coarsest grid, and always shows lower max-norm error without limiting.
The S wave is somewhat less well-resolved, and the difference between
the two curves is smaller, with the max-norm errors both with and
without limiters roughly equal on the coarsest grid.  Finally, the
slow P wave starts out poorly resolved, and using a limiter produces
lower max-norm error on all but the finest grid.  The 1-norm error is
lower with the limiter included for all cases.  For further discussion
of the benefits and drawbacks of using limiters, see
LeVeque~\cite{rjl:fvm-book}.  There have also been efforts to produce
limiters that are compatible with higher-order methods; see for
example \v{C}ada and Torrilhon~\cite{cada-torrilhon:3rdorder-ldlr},
Liu and Tadmor~\cite{liu-tadmor:3rdorder}, or the recent review by
Kemm~\cite{kemm:limiters}.

\subsection{Plane wave convergence results --- viscous}

With viscosity included, numerical solution of the equations of
poroelasticity becomes substantially more challenging.  The chief
difficulty is that the dissipation term has its own associated time
scales, independent of the computational grid.  Since an appropriate
time step for the wave propagation part of the system is proportional
to the grid size --- preferably with a CFL number near 1 --- for large
enough grid cell sizes the dissipation term is stiff relative to the
wave propagation term.  While we maintain stability by solving the
dissipation term exactly, because the wave propagation and dissipation
parts of the system do not commute, we can still heuristically expect
problems in our operator splitting scheme if the time step is much
longer than the characteristic time scale for dissipation.

\begin{figure}
  \begin{center}
    \subfloat[]{\includegraphics[width=0.45\textwidth]{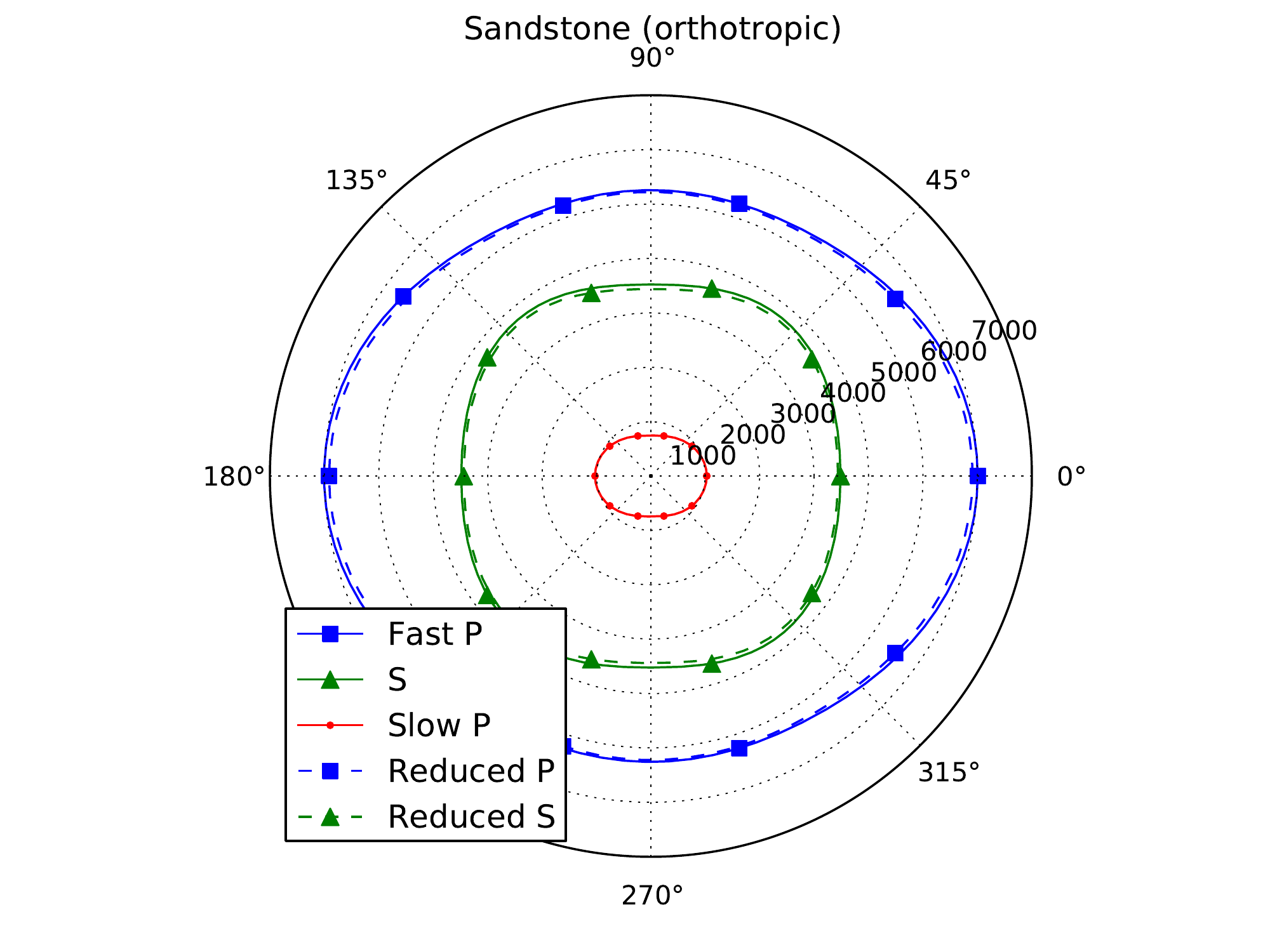}}
    \subfloat[]{\includegraphics[width=0.45\textwidth]{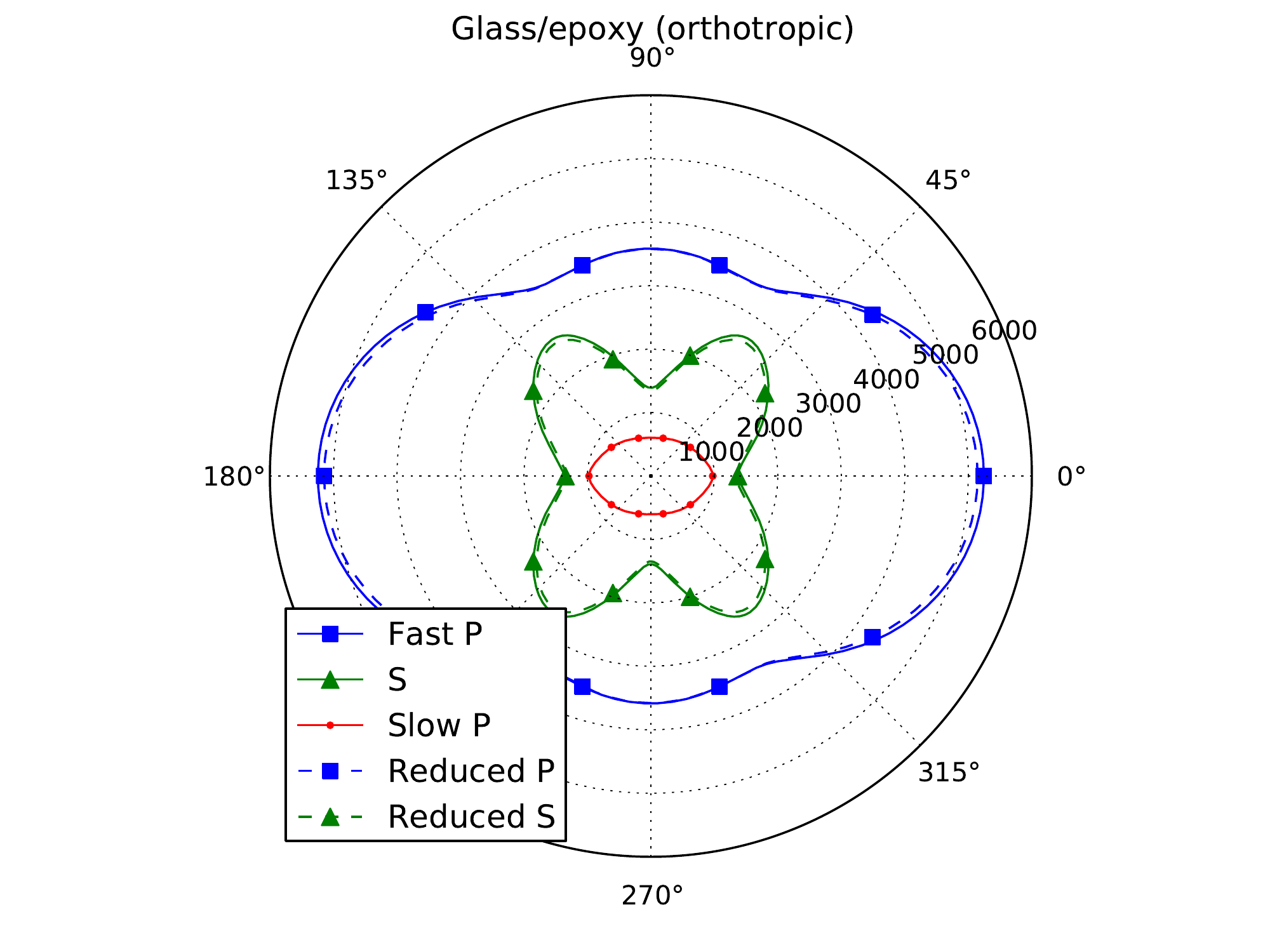}}
    \caption{Wave speeds (meters per second) for full and reduced
      systems for both materials used with viscosity included, as a
      function of propagation angle.  Zero propagation angle
      corresponds to the material principal $1$-axis.  The strict
      subcharacteristic condition \eqref{eq:subcharacteristic} is
      satisfied.
    \label{fig:subcharacteristic}}
  \end{center}
\end{figure}

Revisiting the subcharacteristic condition of section
\ref{sec:stiff-subcharacteristic}, Figure \ref{fig:subcharacteristic}
shows the wave speeds for the full and reduced systems as a function
of propagation direction relative to the principal axes.  The strict
subcharacteristic condition \eqref{eq:subcharacteristic} is satisfied
for all the materials we examine, although the reduced P wave speed
nearly reaches the fast P wave speed for the glass/epoxy material.
In fact, an even stricter condition is satisfied: the wave speeds
interleave, with exactly one wave of the reduced system between each
consecutive pair of waves of the full system.  Based on the discussion
of section \ref{sec:stiff-subcharacteristic}, this suggests that we
will not see spurious solutions or incorrect wave speeds from our
numerical solution.

To better explore the behavior of our numerical method in the presence
of dissipation, we ran a series of numerical tests in the same
sandstone medium as the inviscid test cases, against plane wave
solutions at frequencies ranging from 10 Hz to 20 kHz.  (The maximum
frequency for low-frequency Biot theory to be valid in this medium is
roughly 25 kHz.)  For all cases, the material 1-axis was aligned with
the global $x$-axis, and waves were set to propagate in the positive
$x$ direction.  For the fast P and S waves, we chose the domain size
to be two damped wavelengths of the wave in question; since the slow P
wave has a characteristic decay length (the distance over which the
wave amplitude decreases by a factor of $e$) that is typically a fifth
or less of its wavelength, we chose the domain size for the slow P
wave cases to be twice the characteristic decay length instead.  We
chose the domain sizes in this way so that the number of grid cells
per wavelength, or per decay length, would be constant across all
frequencies; this keeps the discretization error contributed from the
wave propagation part of the system roughly constant for each grid
size across all frequencies, helping to isolate the error caused by
operator splitting.  The total simulation time for the fast P and S
wave cases was 1.25 cycles of the wave (we chose a non-integer number
of cycles to avoid any possible spoofing where an unchanged solution
might appear correct), while for the slow P wave cases it was 1.25
times the time for a fast P wave traveling in the $x$ direction to
cross the domain.  For each combination of wave family and frequency,
we computed solutions on grids of size $100 \times 100$, $200 \times
200$, $400 \times 400$, and $800 \times 800$, using both Godunov and
Strang splitting.

\begin{figure}
  \begin{center}
    \includegraphics[width=0.8\textwidth]{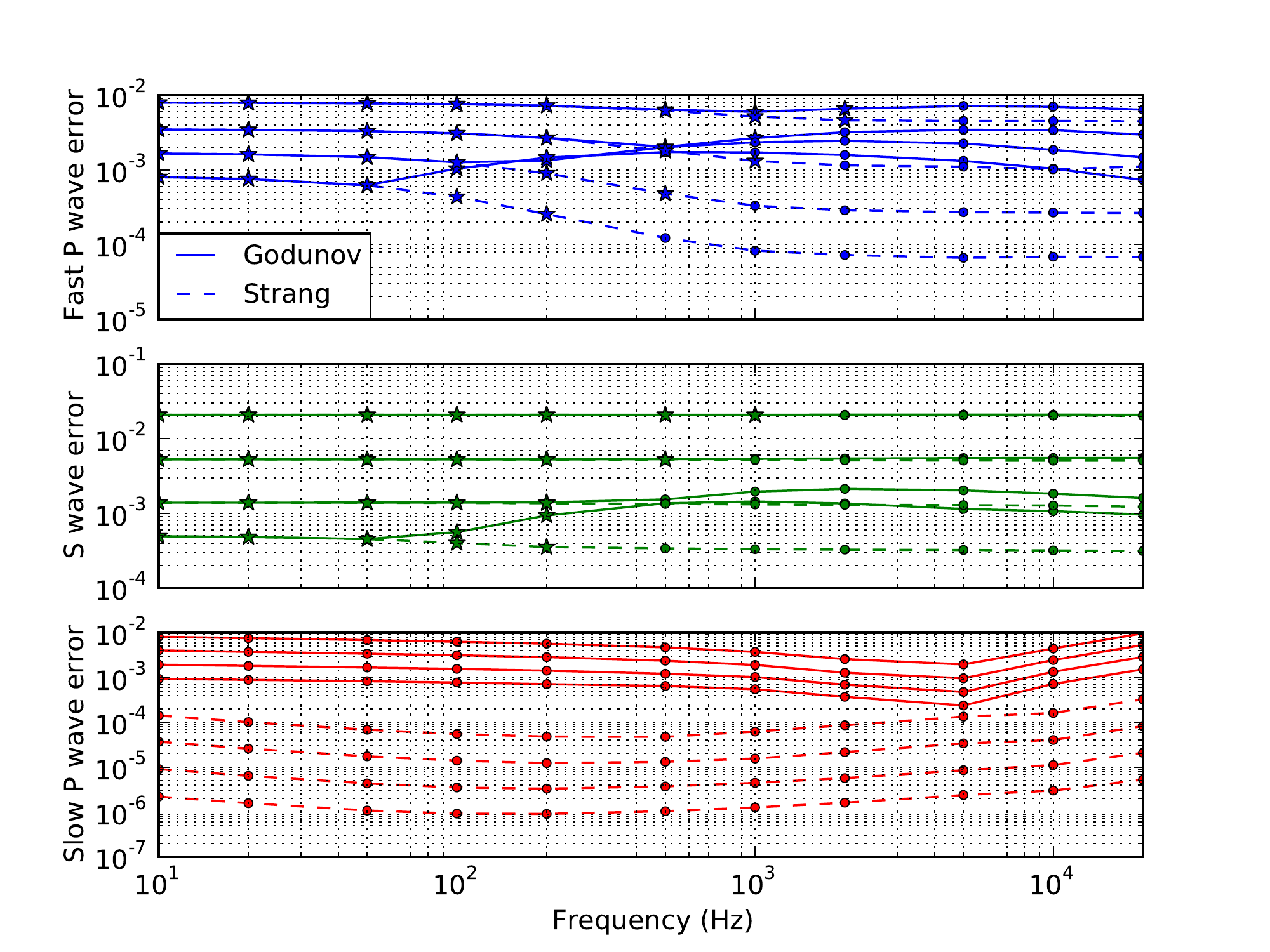}
    \caption{Normalized energy max-norm error for operator splitting
      tests.  For each splitting method (Godunov or Strang) within
      each subplot, the curves fall in order of increasing grid
      fineness, with each curve from a grid twice as fine as the curve
      above.  The top curve is on a $100 \times 100$ grid, and the
      bottom is on an $800 \times 800$.  The circles indicate where
      the time step was less than the characteristic time scale $\tau_d$
      for dissipation in the $x$ direction; the stars indicate where
      the time step was greater.
      \label{fig:opsplit-results}}
  \end{center}
\end{figure}

Figure \ref{fig:opsplit-results} shows the results of these tests in
the same normalized energy max-norm used for the inviscid cases.
There is a pronounced qualitative difference in convergence behavior
depending on frequency.  At low frequencies, corresponding to large
grid cell sizes and long time steps, both splitting methods show
first-order convergence.  Starting at a step length of roughly 5-10
times the characteristic time $\tau_d$ for dissipation in the $x$
direction (the time over which the $x$ velocity of the fluid relative
to the matrix decreases by a factor of $e$), the two methods begin
behaving differently, with the Godunov splitting error increasing
abruptly while the Strang splitting error sweeps smoothly down to
second-order convergence.  This effect is most visible for the fast P
wave, but can also be seen slightly in the $400 \times 400$ grid and
strongly in the $800 \times 800$ grid for the S wave.  Because of the
choice of domain size, the slow P test cases always had a time step
below the characteristic dissipation time; they display consistent
first-order convergence with Godunov splitting and second-order with
Strang.  The qualitative shift in behavior with time step length can
be understood by noting that for a time step much longer than the
characteristic time, the solution operator $\exp(\bD\Delta t)$ of the
dissipative part of the system is essentially a projection operator
that sets the fluid relative velocity to zero and transfers all of the
fluid relative momentum into the bulk motion of the medium.  The
effect of this projection operator is essentially the same whether it
is applied once per time step after solution of the wave propagation
part of the system (Godunov splitting), or twice, both before and
after wave propagation (Strang splitting).

With these results available to inform our choices, we conducted a set
of convergence studies similar to the inviscid cases in the previous
subsection.  Because of the qualitative difference in convergence
behavior for different time step regimes, we performed convergence
studies both at a point in the high-frequency convergence regime of
Figure \ref{fig:opsplit-results} (10 kHz), and at a low-frequency
point (10 Hz).  Since the slow P wave decays extremely rapidly in the
presence of viscosity --- typically by a factor of 10 to 100 or more
per wavelength in the valid frequency range for Biot theory --- we
chose a different domain size for the viscous slow P test cases.
For the high-frequency cases, we chose square domains with side length
1.2 meters for the fast P and S waves --- roughly two wavelengths of
the fast P wave at 10 kHz --- and 5 centimeters for the slow P wave ---
roughly 2-5 times the characteristic decay length for this wave,
depending on propagation direction in the orthotropic medium.  For the
low-frequency cases, the domain size was 1200 meters for the fast P
and S waves, and 1 meter for the slow P wave --- again, 2-4 times the
characteristic decay length, which is far shorter than the
wavelength.  This huge disparity in domain sizes is somewhat
troublesome, but it is not clear whether simulation results for a slow
P wave on a domain of the size used for the fast P and S waves would
be meaningful, since the solution decays over such a short distance.
For practical problems, this would be an excellent opportunity for
adaptive mesh refinement, to generate fine grids where and when slow P
waves appear, then coarsen the grid again after they dissipate.

We ran all of the viscous test cases to essentially the same final
times as for the frequency sweep of Figure \ref{fig:opsplit-results}.
For the $10\,\text{kHz}$ runs, this was $125\,\text{$\mu$s}$ for the fast P and S
waves (1.25 periods of the wave), and $10.4\,\text{$\mu$s}$ for the slow P
wave (1.25 times the time for a fast P wave to cross the domain),
while for the $10\,\text{Hz}$ runs it was $0.125\,\text{s}$ for the fast P and S
waves, and $208\,\text{$\mu$s}$ for the slow P wave.  Strang splitting was
used for all cases.  All other aspects of the solution, including the
sets of wave propagation and material principal directions
$\theta_{\text{wave}}$ and $\theta_{\text{mat}}$ as well as the method
of setting the boundary condtions, were the same as for the inviscid
test cases, and we measured the solution error using the same set of
norms.

Table \ref{tab:conv-summary-visc-hifreq} summarizes convergence for
the high-frequency viscous cases.
We again see consistent second-order convergence in all norms for the
fast P and S waves at high frequency, and a similar amount of
dependence of error on wave propagation direction.  For the slow P
wave at high frequency, however, results are substantially different
from the inviscid cases.  We see second-order convergence, since the
solution is well-resolved on the grids used here, but there is now a
factor of several hundred difference between the maximum and minimum
error at a single grid size.  Close examination of the error for
individual cases shows that it is primarily a function of the offset
$\theta_{\text{wave}} - \theta_{\text{mat}}$; for a fixed value of
$\theta_{\text{wave}} - \theta_{\text{mat}}$ and a fixed grid size,
the error is similar across all values of $\theta_{\text{mat}}$.  This
indicates that the large variation in error is a effect of the
alignment of the wavefront relative to the principal material axes,
rather than a grid alignment effect.  The likely cause is the
substantial difference in the characteristic decay times between the 1
and 3 axes of the material --- the decay time in the 1 direction is
5.95 microseconds, while in the 3 direction it is 1.82 microseconds.
This large variation in decay time causes a large variation in the
operator splitting error.  In addition, the characteristic decay
length is substantially shorter in the 3 direction, causing the
solution magnitude to be larger at the ``upstream'' (opposite the
propagation direction) edge of the domain relative to the value at the
center against which the error is normalized; the larger solution
magnitude naturally results in a larger error.

\begin{table}
\caption{Summary of convergence results for viscous, high-frequency
  test cases}
\label{tab:conv-summary-visc-hifreq}
\begin{center}
\begin{tabular}{cccccccc}
  & & \multicolumn{3}{c}{Convergence rate} & &
  \multicolumn{2}{c}{Error on $800 \times 800$ grid}\\
  \cline{3-5} \cline{7-8}
  Wave family & Error norm & Best & Worst & Mean & Worst $R^2$ value &
  Best & Worst \\
  \hline
  \multirow{3}{*}{Fast P} & 1-norm & 2.03 & 2.01 & 2.01 & 0.99996 &
  $2.61 \times 10^{-5}$ & $7.31 \times 10^{-5}$ \\
  & 2-norm & 2.03 & 2.01 & 2.01 & 0.99996 & $3.16 \times 10^{-5}$ &
  $8.79 \times 10^{-5}$ \\
  & Max-norm & 2.05 & 2.00 & 2.02 & 0.99887 & $6.44 \times 10^{-5}$ &
  $1.77 \times 10^{-4}$ \\
  \hline
  \multirow{3}{*}{S} & 1-norm & 2.01 & 2.01 & 2.01 & 1.00000 &
  $1.38 \times 10^{-4}$ & $3.39 \times 10^{-4}$ \\
  & 2-norm & 2.01 & 2.00 & 2.01 & 1.00000 & $1.58 \times 10^{-4}$ &
  $3.98 \times 10^{-4}$ \\
  & Max-norm & 2.03 & 1.99 & 2.00 & 0.99978 & $2.99 \times 10^{-4}$ &
  $7.59 \times 10^{-4}$ \\
  \hline
  \multirow{3}{*}{Slow P} & 1-norm & 2.01 & 2.00 & 2.01 & 1.00000 &
  $2.12 \times 10^{-6}$ & $1.74 \times 10^{-4}$ \\
  & 2-norm & 2.02 & 2.00 & 2.01 & 1.00000 & $2.49 \times 10^{-6}$ &
  $3.26 \times 10^{-4}$ \\
  & Max-norm & 2.02 & 1.96 & 2.00 & 0.99988 & $6.53 \times 10^{-6}$ &
  $2.25 \times 10^{-3}$
\end{tabular}
\end{center}
\end{table}

Table \ref{tab:conv-summary-visc-lofreq}
shows the results for the
low-frequency viscous test cases.  The fast P and S waves again show
only a weak dependence of error on wave propagation and principal
material direction, but their convergence rates are substantially
degraded, just as in the low-frequency range of Figure
\ref{fig:opsplit-results}.  Convergence of the fast P wave is roughly
first-order in the worst case, due to the long time step relative to
the characteristic decay times.  Surprisingly, results for the S wave
are only slightly better in the worst case, likely due to the much
shorter characteristic decay time in the 3 direction.
If the grid were further refined, we would
presumably reach a second-order convergence regime as the time step
approached the characteristic decay time, but in order to have a time
step similar to the shortest decay time at a CFL number of 0.9 we
would need a grid cell size of roughly $13\,\text{cm}$ --- resulting in
a $90000 \times 90000$ cell grid on the $1200\,\text{m}$ square domain if
this size is uniform!  The low-frequency slow P wave cases, by
contrast, show the same strong dependence of error on the alignment of
the wave direction with the principal material direction as the
high-frequency cases, but because the domain size was much smaller and
the time step much shorter, we observe consistent second-order
convergence.

\begin{table}
\caption{Summary of convergence results for viscous, low-frequency
  test cases}
\label{tab:conv-summary-visc-lofreq}
\begin{center}
\begin{tabular}{cccccccc}
  & & \multicolumn{3}{c}{Convergence rate} & &
  \multicolumn{2}{c}{Error on $800 \times 800$ grid}\\
  \cline{3-5} \cline{7-8}
  Wave family & Error norm & Best & Worst & Mean & Worst $R^2$ value &
  Best & Worst \\
  \hline
  \multirow{3}{*}{Fast P} & 1-norm & 1.57 & 1.10 & 1.26 & 0.99239 &
  $1.95 \times 10^{-4}$ & $4.24 \times 10^{-4}$ \\
  & 2-norm & 1.60 & 1.10 & 1.28 & 0.99206 & $2.15 \times 10^{-4}$ &
  $4.77 \times 10^{-4}$ \\
  & Max-norm & 1.61 & 1.10 & 1.32 & 0.99014 & $3.90 \times 10^{-4}$ &
  $1.04 \times 10^{-3}$ \\
  \hline
  \multirow{3}{*}{S} & 1-norm & 1.85 & 1.37 & 1.58 & 0.99241 &
  $4.37 \times 10^{-4}$ & $7.30 \times 10^{-4}$ \\
  & 2-norm & 1.85 & 1.37 & 1.59 & 0.99274 & $5.01 \times 10^{-4}$ &
  $8.30 \times 10^{-4}$ \\
  & Max-norm & 1.83 & 1.36 & 1.58 & 0.98917 & $9.63 \times 10^{-4}$ &
  $1.92 \times 10^{-3}$ \\
  \hline
  \multirow{3}{*}{Slow P} & 1-norm & 2.00 & 1.94 & 1.98 & 0.99985 &
  $1.04 \times 10^{-6}$ & $4.88 \times 10^{-5}$ \\
  & 2-norm & 2.03 & 1.93 & 1.99 & 0.99982 & $1.08 \times 10^{-6}$ &
  $7.40 \times 10^{-5}$ \\
  & Max-norm & 2.09 & 1.84 & 1.98 & 0.99795 & $2.40 \times 10^{-6}$ &
  $4.59 \times 10^{-4}$
\end{tabular}
\end{center}
\end{table}

As a final comment on convergence, we note that while the rate of
reduction of error with decreasing mesh size is poorer for the
low-frequency fast P and S wave cases, this may not be a problem in
practice unless very high accuracy is desired.  For both frequency
ranges investigated, the relative error in all three norms on the $800
\times 800$ grid never exceeds $1.04 \times 10^{-3}$ for the fast P
wave, or $1.92 \times 10^{-3}$ for the S wave --- despite the numerical
difficulties encountered, the waves are still well-resolved.

\subsection{Single-material point source results}

With the accuracy of our method characterized for simple plane wave
solutions, we are able to move on to more interesting problems.
First, we compare against the results of de la Puente et
al.~\cite{delapuente-dumbser-kaser-igel:poro-dg} and
Carcione~\cite{carcione:poro-aniso-1996} for a point source in a
uniform orthotropic medium.  We used the orthotropic sandstone and
glass/epoxy media described in Table \ref{tab:matprops}; in both
cases, the material 1 axis coincided with the $x$ axis.  All test
cases started with initial condition $\bQ(x,z,0) = 0$, and the domain
was excited by a point source with a Ricker wavelet profile having peak
frequency $f_{\text{src}} = 3730 \, \text{Hz}$ for sandstone and $3135
\, \text{Hz}$ for glass/epoxy, acting with peak intensity $+1\,
\text{Pa}\cdot\text{m}^2/\text{s}$ on the vertical normal stress
$\sigma_{zz}$ and $-1\, \text{Pa}\cdot\text{m}^2/\text{s}$ on fluid
pressure.  The peak of the wavelet occurred $0.4 \,\text{ms}$ after
the start of the simulation.  For each case, we used a square domain
$18.7\,\text{m}$ on a side; the point source was placed at the center of
the domain.  The simulation time span was $1.56 \, \text{ms}$ for the
sandstone medium, and $1.80 \, \text{ms}$ for glass/epoxy.  We
calculated results both with and without viscosity present.

\begin{figure}
  \begin{center}
    \subfloat[]{\label{fig:sandstone-invisc-vx}
      \includegraphics[width=0.24\textwidth, trim=0.5in 0 1in 0, clip]{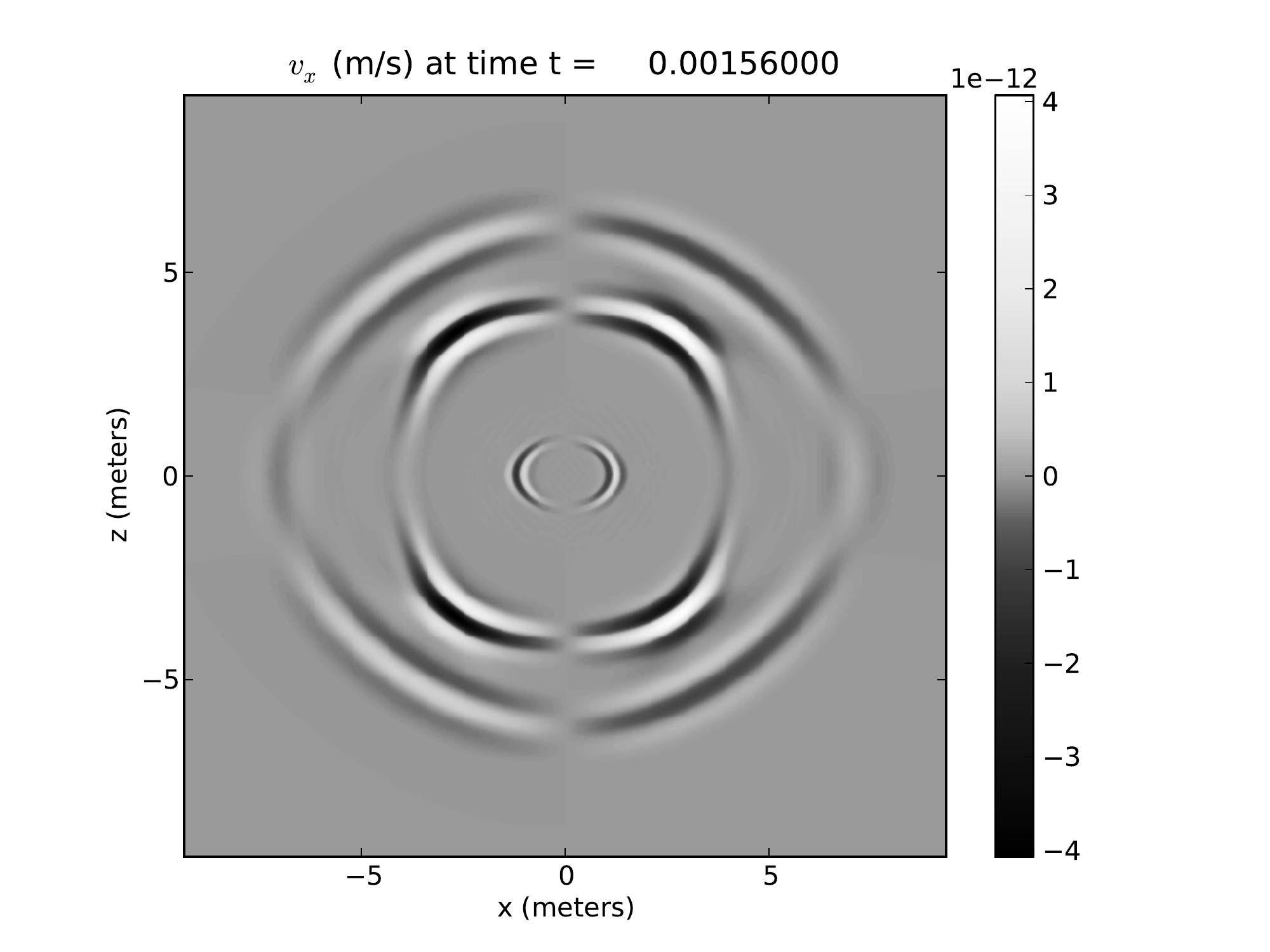}}
    \subfloat[]{\label{fig:sandstone-visc-vx}
      \includegraphics[width=0.24\textwidth, trim=0.5in 0 1in 0, clip]{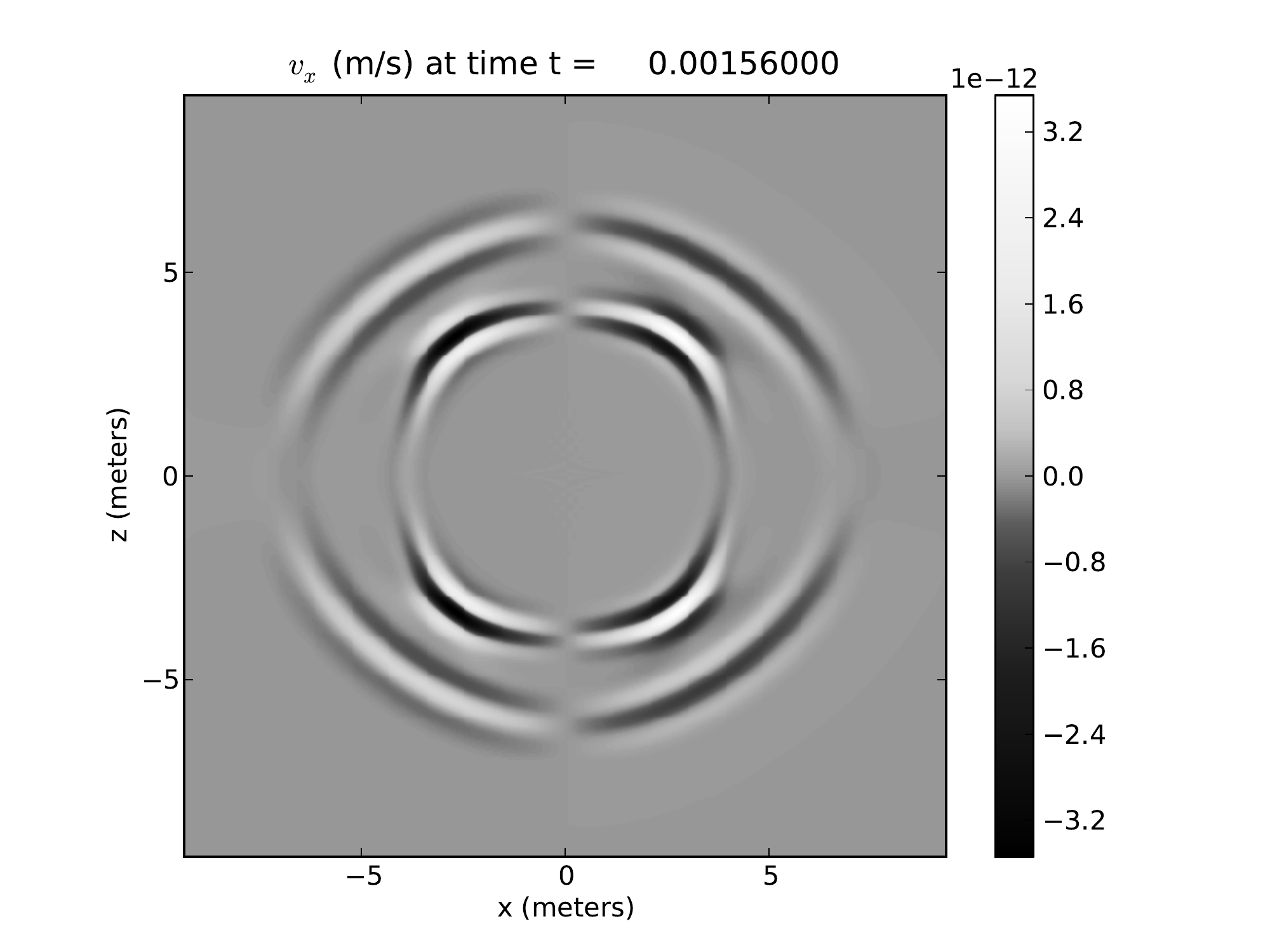}}
    \subfloat[]{\label{fig:sandstone-invisc-vz}
      \includegraphics[width=0.24\textwidth, trim=0.5in 0 1in 0, clip]{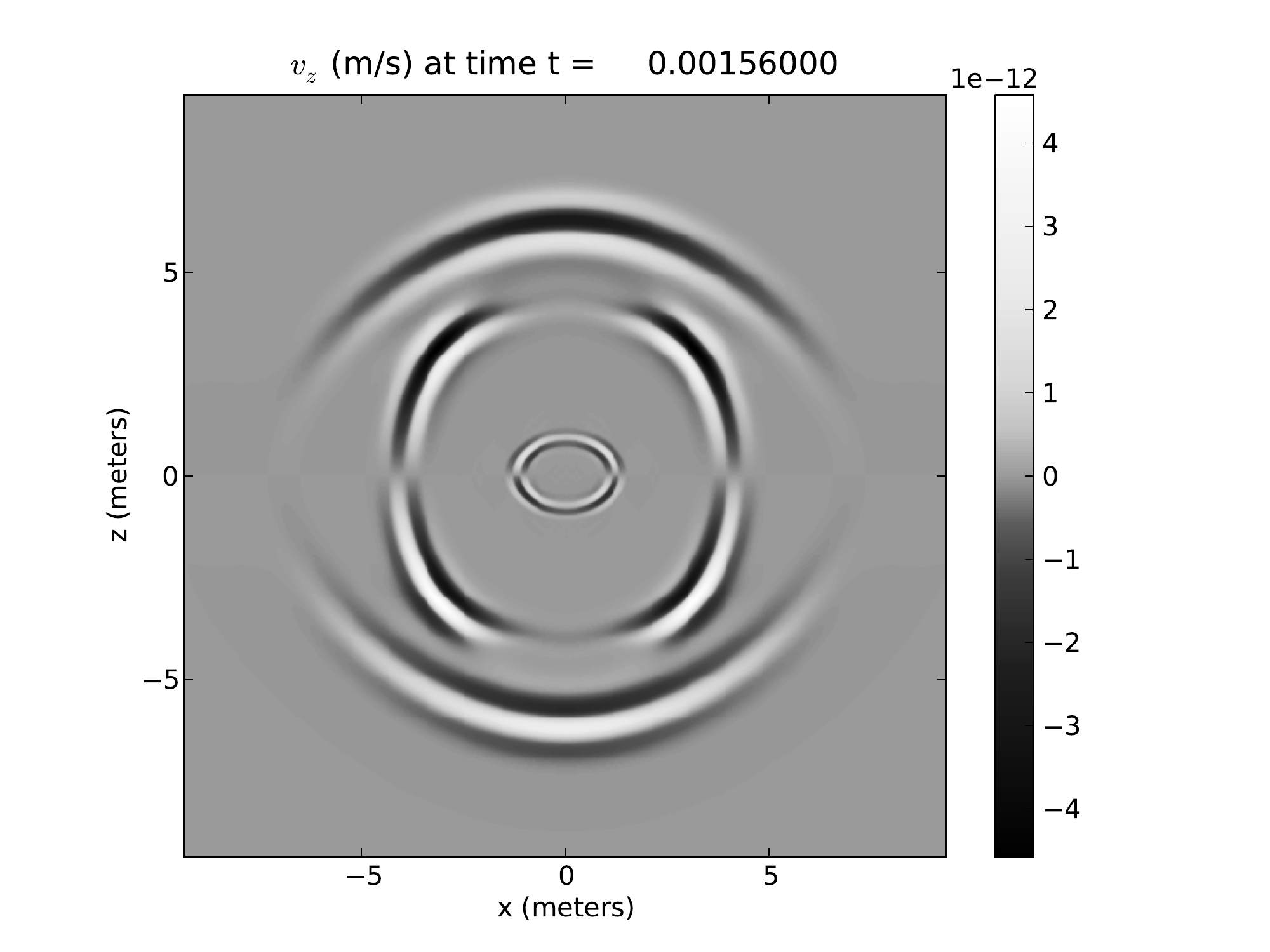}}
    \subfloat[]{\label{fig:sandstone-visc-vz}
      \includegraphics[width=0.24\textwidth, trim=0.5in 0 1in 0, clip]{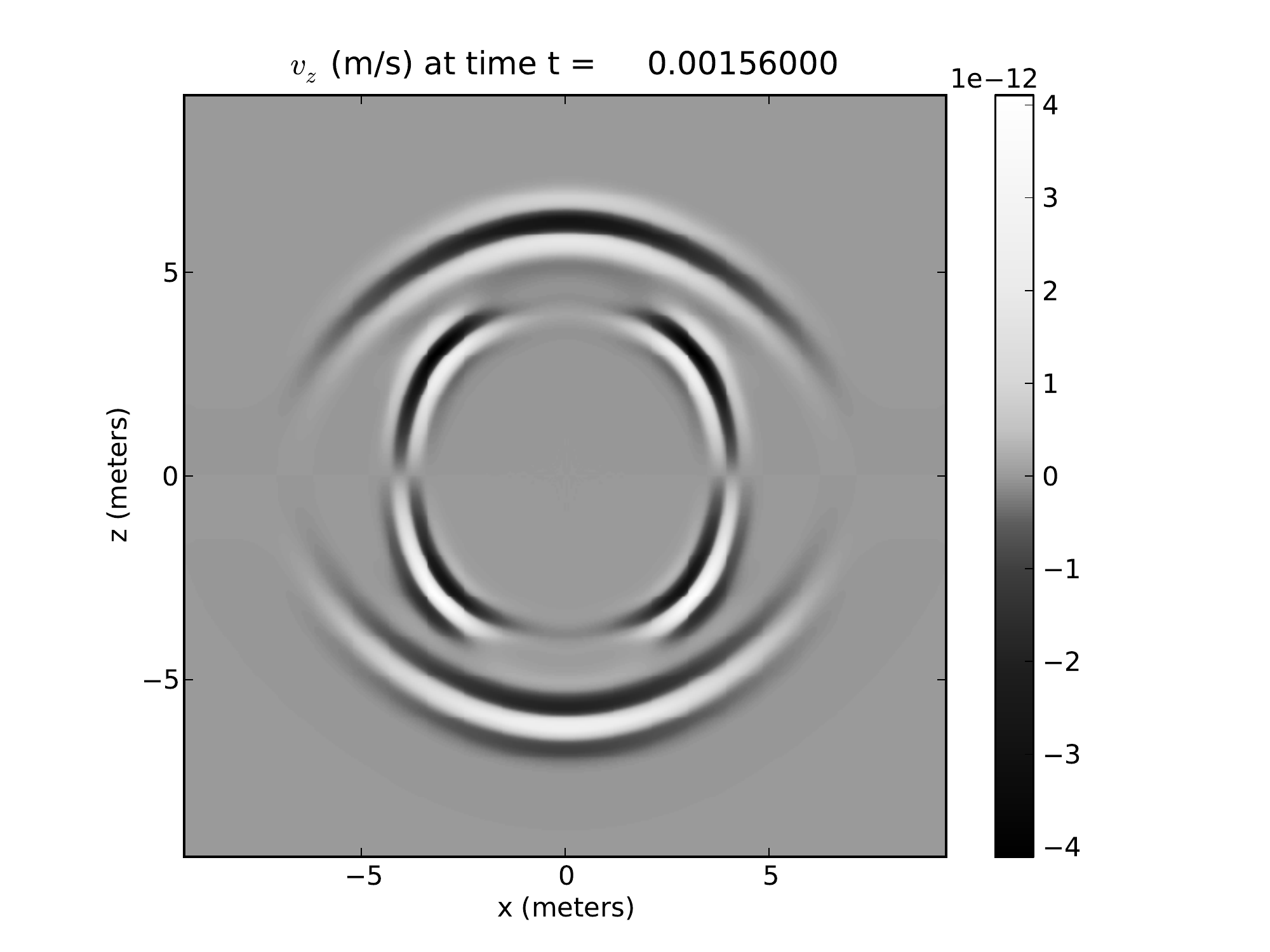}}\\
    \subfloat[]{\label{fig:gfrp-invisc-vx}
      \includegraphics[width=0.24\textwidth, trim=0.5in 0 1in 0, clip]{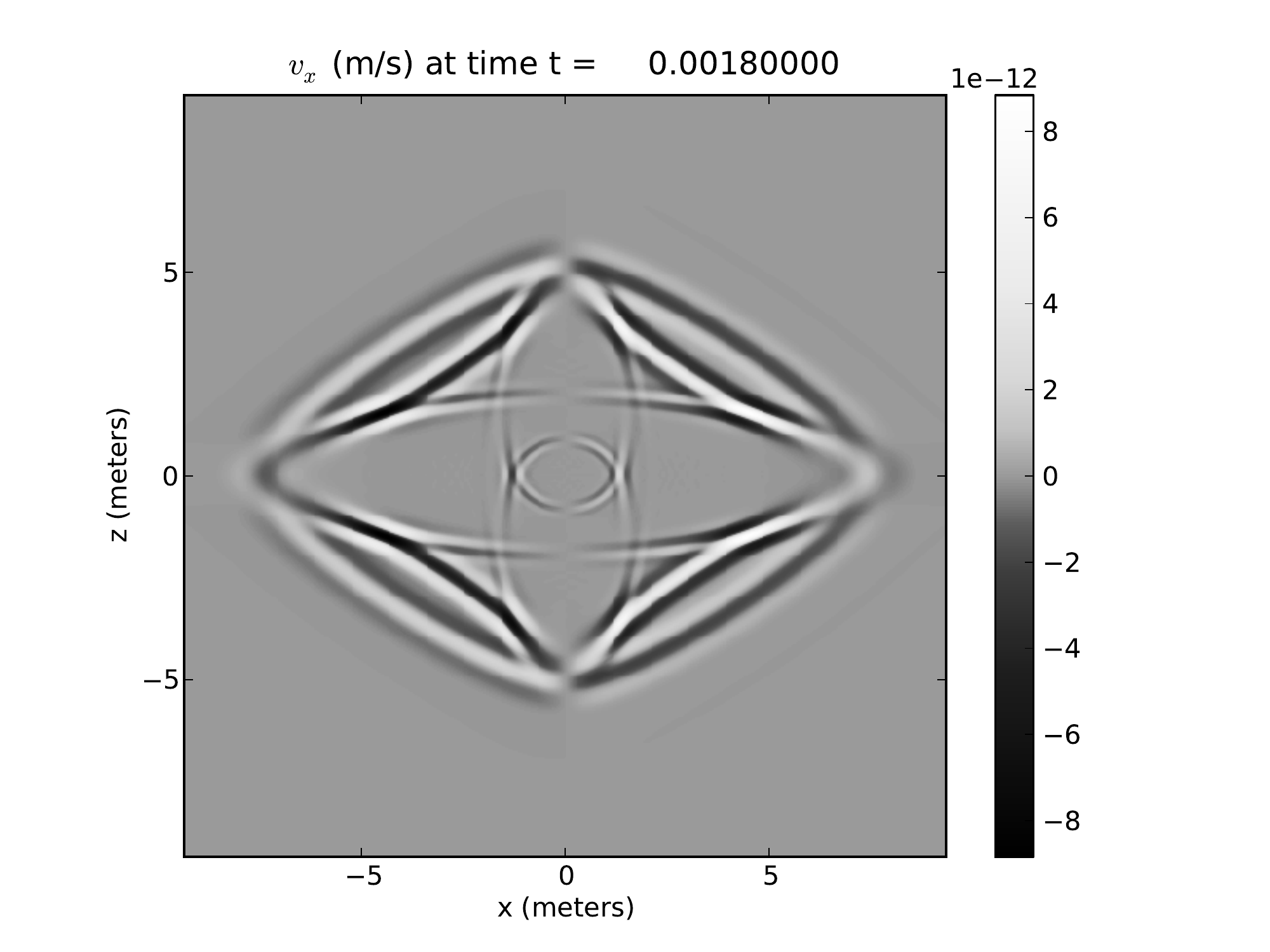}}
    \subfloat[]{\label{fig:gfrp-visc-vx}
      \includegraphics[width=0.24\textwidth, trim=0.5in 0 1in 0, clip]{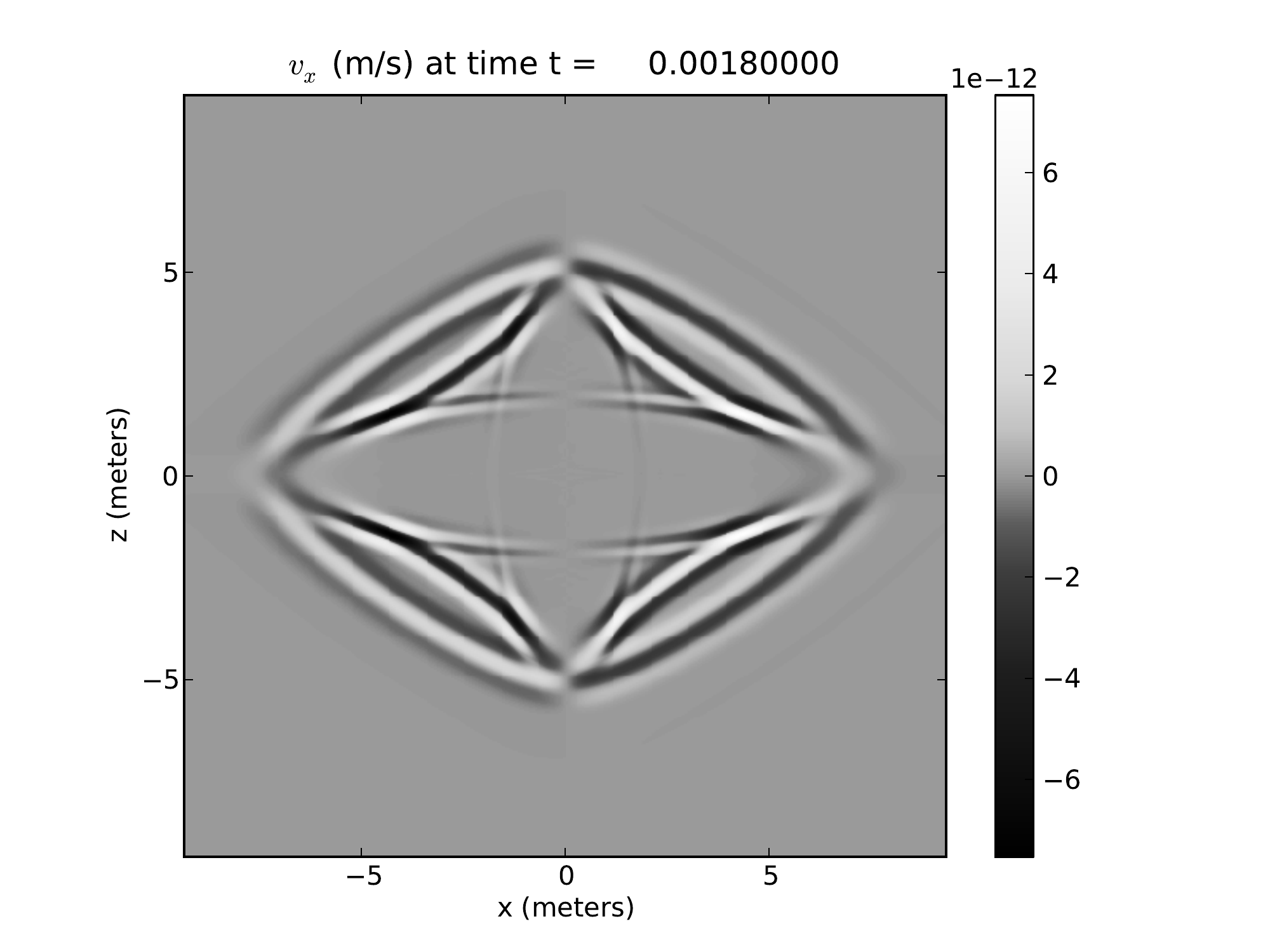}}
    \subfloat[]{\label{fig:gfrp-invisc-vz}
      \includegraphics[width=0.24\textwidth, trim=0.5in 0 1in 0, clip]{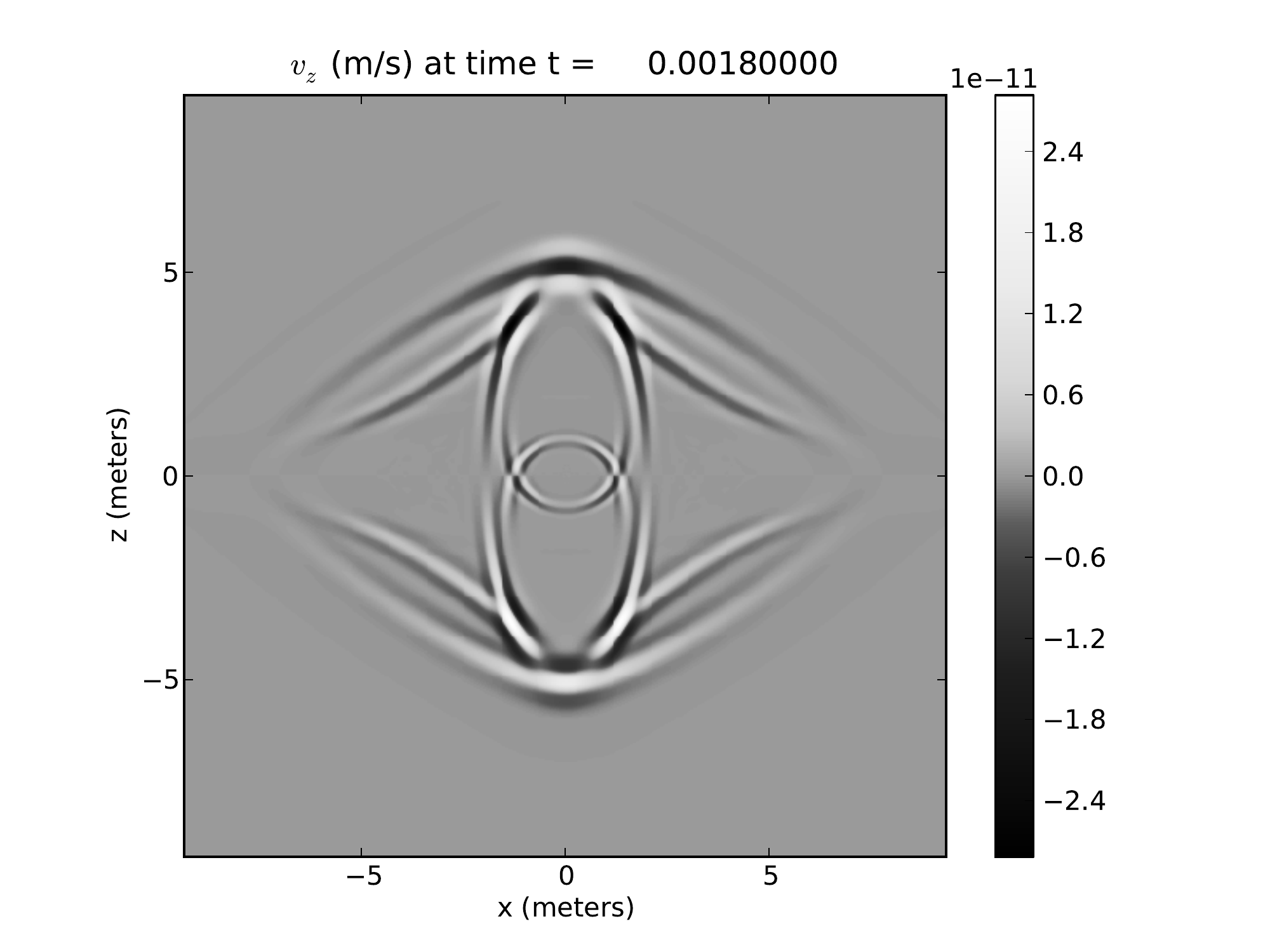}}
    \subfloat[]{\label{fig:gfrp-visc-vz}
      \includegraphics[width=0.24\textwidth, trim=0.5in 0 1in 0, clip]{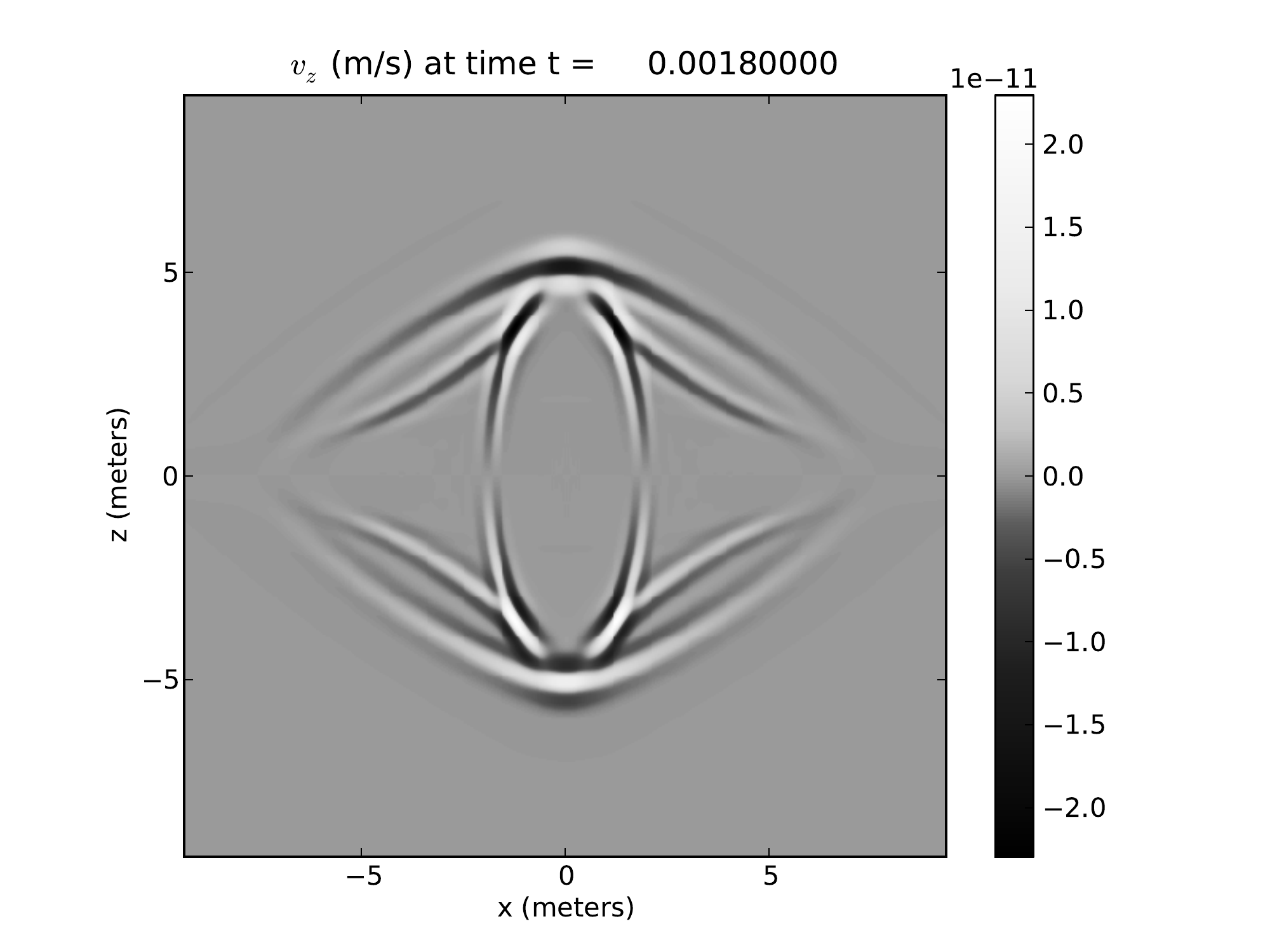}}
    \caption{Results for point source test cases.  Top: sandstone;
      bottom: glass/epoxy; left
      to right: $v_x$ without viscosity, $v_x$ with viscosity, $v_z$
      without viscosity, $v_z$ with viscosity. \label{fig:results-pointsrc}}
  \end{center}
\end{figure}

We carried out all of these simulations on a uniform $501 \times 501$
cell grid.  The odd number of grid cells in each axis allowed us to
apply the point source to a single grid cell.  A grid resolution this
high was necessary to resolve the slow P wave well; as in the plane
wave test cases, the fast P and S waves were well-resolved on
substantially coarser grids.  The point source was implemented
numerically as part of the source step in the operator splitting
scheme; since the point source acts on stress variables, and the
viscous dissipation acts on velocity variables, it does not matter
which is applied first.  The CFL number for all simulations was 0.9,
resulting in 279 timesteps being taken for the sandstone case and 282
for the glass/epoxy.  We used the monotonized centered (MC) limiter
here --- even though most of the solution is smooth, without a limiter
the lack of smoothness at the source point produces substantial
spurious oscillations in the solution.  Figure
\ref{fig:results-pointsrc} shows the results of these simulations.
These figures correspond to Figures 6 and 7 of de la Puente et
al.~\cite{delapuente-dumbser-kaser-igel:poro-dg}, or Figures 5 and 7
of Carcione~\cite{carcione:poro-aniso-1996}, and are in agreement with
them.

\subsection{Heterogeneous domain results}

We conclude our results here with a wave reflection and
interconversion problem that demonstrates the ability of our code to
model material interfaces, as well as the benefits of using adaptive
mesh refinement --- a feature of the \amrclaw{} variant of \clawpack{} that
we have not used thus far.

Our final test case is a large-scale inviscid problem with a bed of
isotropic shale overlying isotropic sandstone, with material
properties given in Table \ref{tab:matprops}.  The problem domain is
the rectangle $[0, 1500\,\text{m}] \times [0, 1400\,\text{m}]$ in the
$x$-$z$ plane, with the boundary between the materials at $z =
700\,\text{m}$.  There is a point source at $(x,z) = (750\,\text{m},
900\,\text{m})$, again with a Ricker wavelet profile in time, with
peak frequency $50 \, \text{Hz}$.  The source acts on the
$z$-direction normal stress and fluid pressure with peak intensities
$+2.3\times 10^{13}\, \text{Pa}\cdot\text{m}^2/\text{s}$ and
$-2.3\times 10^{13}\, \text{Pa}\cdot\text{m}^2/\text{s}$,
respectively, similarly to the previous test cases in homogeneous
domains.  This source magnitude was chosen to match the magnitude of
the response shown by de la Puente et
al.~\cite{delapuente-dumbser-kaser-igel:poro-dg}.  The peak of the
source was delayed $40 \, \text{ms}$ after the start of the
simulation, and the total duration of the run was $0.5 \, \text{s}$.
In addition to time-snapshots of the solution at particular instants,
we also recorded solution time histories at three ``gauges'', located
at $(x_1,z_1) = (950 \,\text{m}, 750\,\text{m})$, $(x_2,z_2) = (950
\,\text{m}, 650\,\text{m})$, and $(x_3,z_3) = (950 \,\text{m},
500\,\text{m})$.

For our adaptively refined simulation, the coarsest-level
computational grid was $75 \times 70$ cells in size, giving square
cells $20\,\text{m}$ on a side.  We used two additional levels of
refinement on this grid, the first at a factor of $4$, and the second
at further factor of $6$, so that cells on the finest grids were
$0.83\,\text{m}$ on a side.  Our code flagged a cell for mesh refinement
when the energy norm (using the material properties of that cell) of
the difference $\Delta \bQ$ between its state vector and that of any
adjacent cell exceeded $32.5 \,\text{J}^{1/2}/\text{m}^{3/2}$,
with the exception of the rectangle $[700\,\text{m}, 1000\,\text{m}]
\times [450\,\text{m}, 950\,\text{m}]$, where the threshold for
refinement was lowered to $3.25 \,\text{J}^{1/2}/\text{m}^{3/2}$
in order to improve accuracy at the gauges.  These tolerances were
chosen empirically based on the observed magnitude of the waves in
the simulation.
This refinement criterion is a generalization of the typical
\amrclaw{} approach of refining based on the difference between the
solution values in neighboring cells, which has been used successfully
on many problems; \amrclaw{} makes it easy to set alternate
user-specified refinement criteria if desired, and also offers
automatic error estimation via Richardson extrapolation.
Besides refinement based on the solution field, the source location
was also
flagged for refinement to the finest level available whenever the
source intensity was greater than about $10^{-9}$ of peak.  Since the
grid size was even in each direction on all but the coarsest grids, the
source was distributed over the four cells closest to its location
using a bilinear weighting.  We again used the MC limiter for this
problem.

\begin{figure}
  \begin{center}
    \subfloat[]{\label{fig:2mat-snapshot-vz}
      \includegraphics[width=0.4\textwidth]{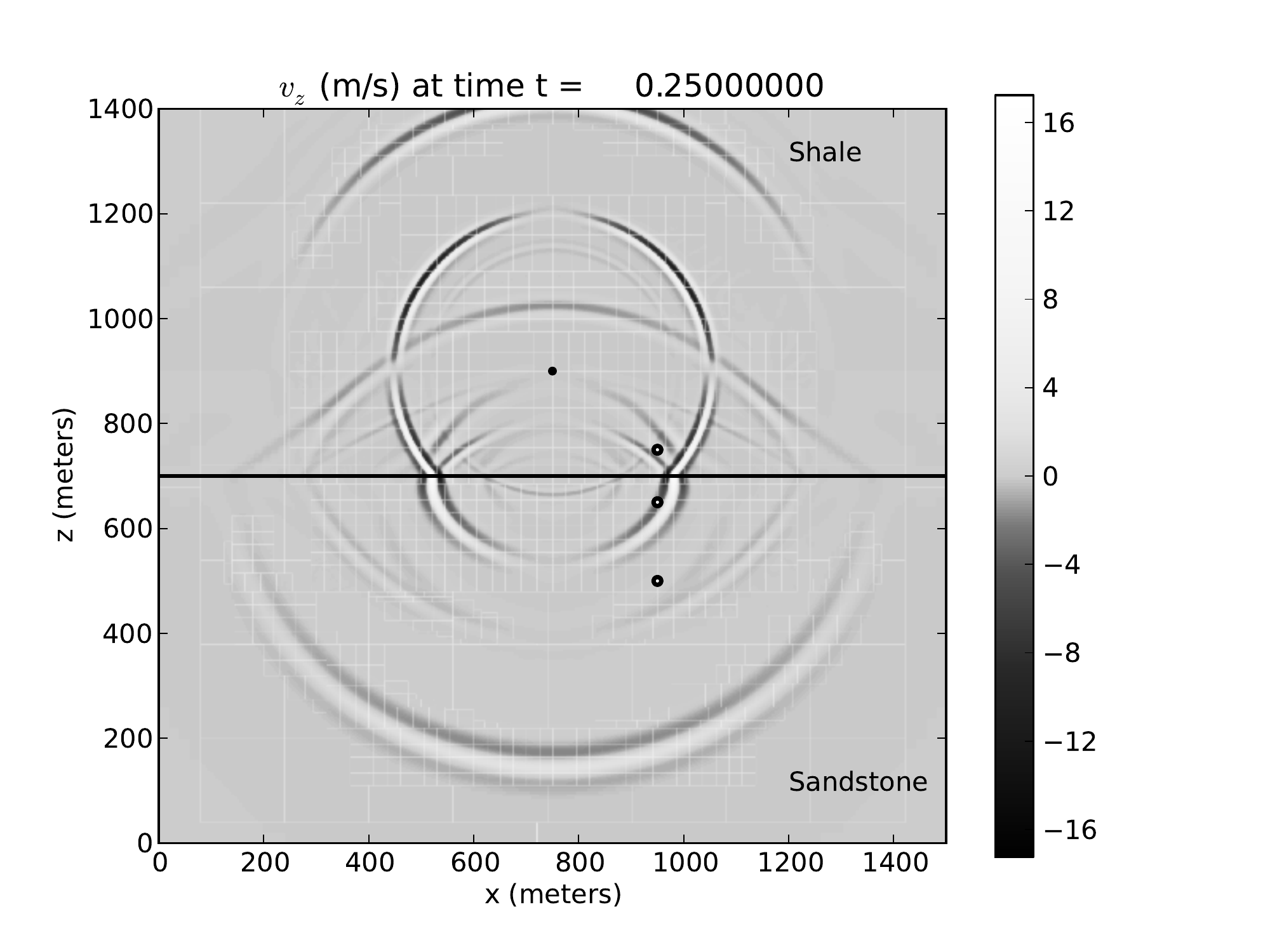}}
    \subfloat[]{\label{fig:2mat-snapshot-grids}
      \includegraphics[width=0.4\textwidth]{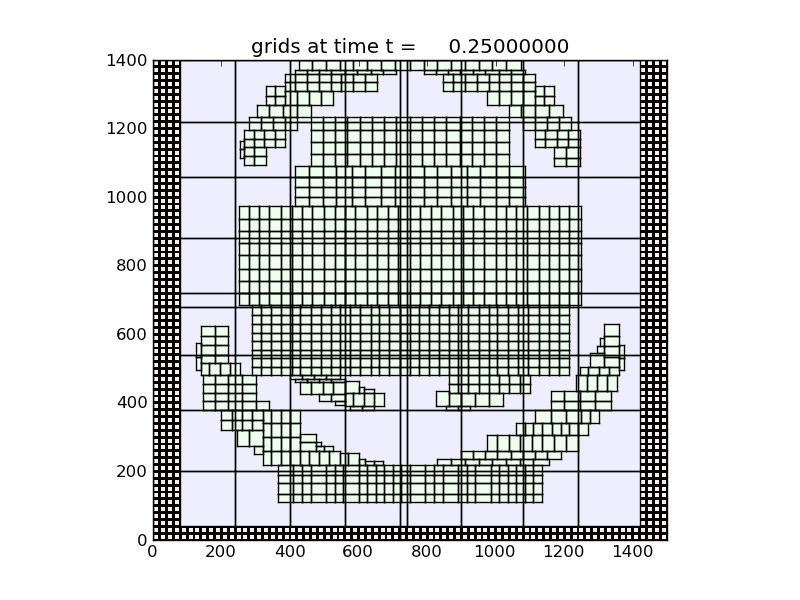}}
    \caption{Left: snapshot of $z$-direction solid velocity 0.25
      seconds after the start of the simulation.  The source location
      is marked with a solid black dot, while the gauges are marked
      with white-centered black dots.  The gauges are numbered from
      top to bottom.  Right: AMR grids at this time point.  Individual
      cells are drawn on the coarsest AMR level, but only grid
      outlines are shown on finer levels. \label{fig:2mat-snapshot}}
  \end{center}
\end{figure}

\begin{figure}
  \begin{center}
    \subfloat[]{\includegraphics[width=0.3\textwidth,
        trim=0.34in 0 0.6in 0, clip]{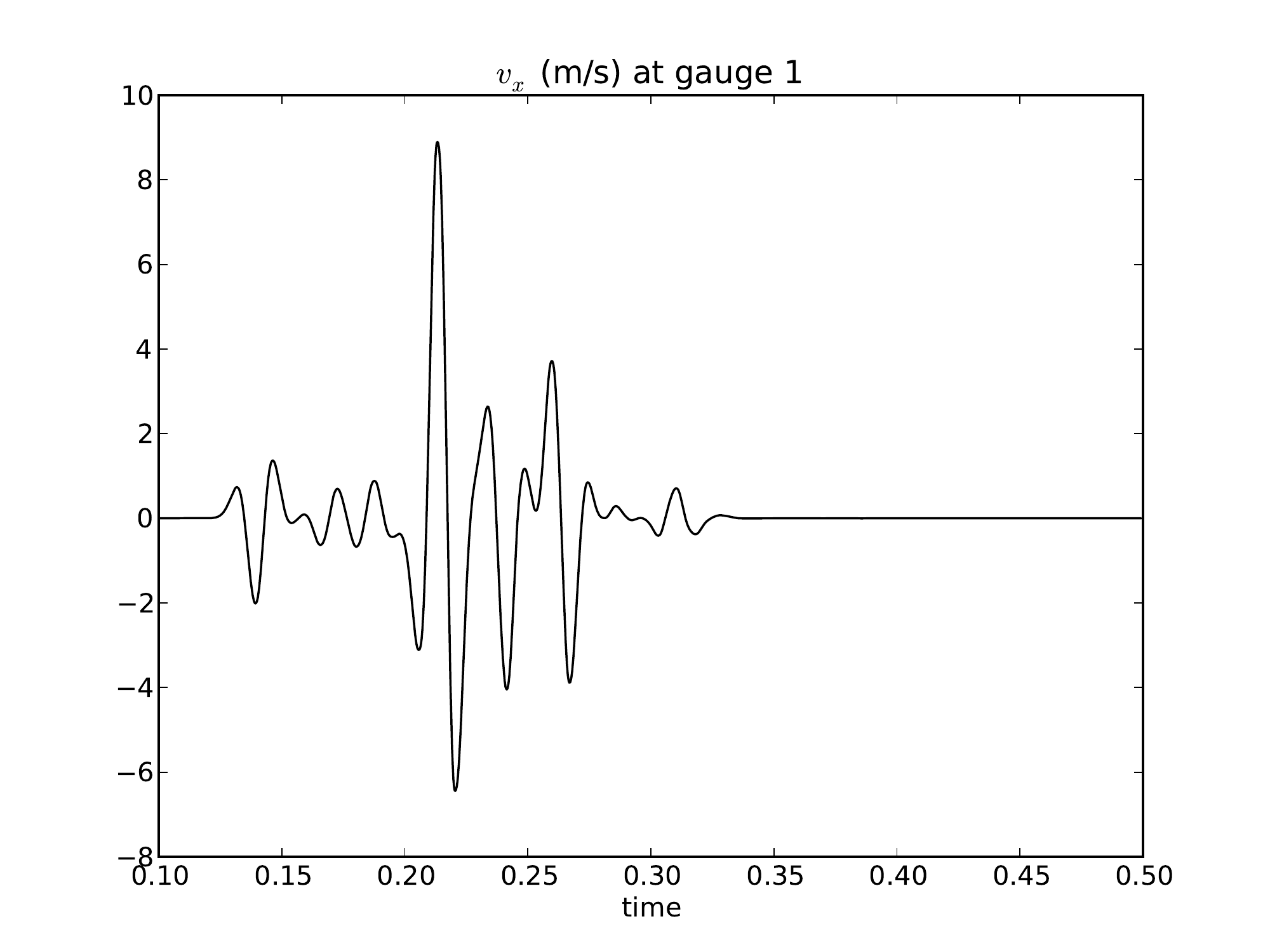}}
    \subfloat[]{\includegraphics[width=0.3\textwidth,
        trim=0.34in 0 0.6in 0, clip]{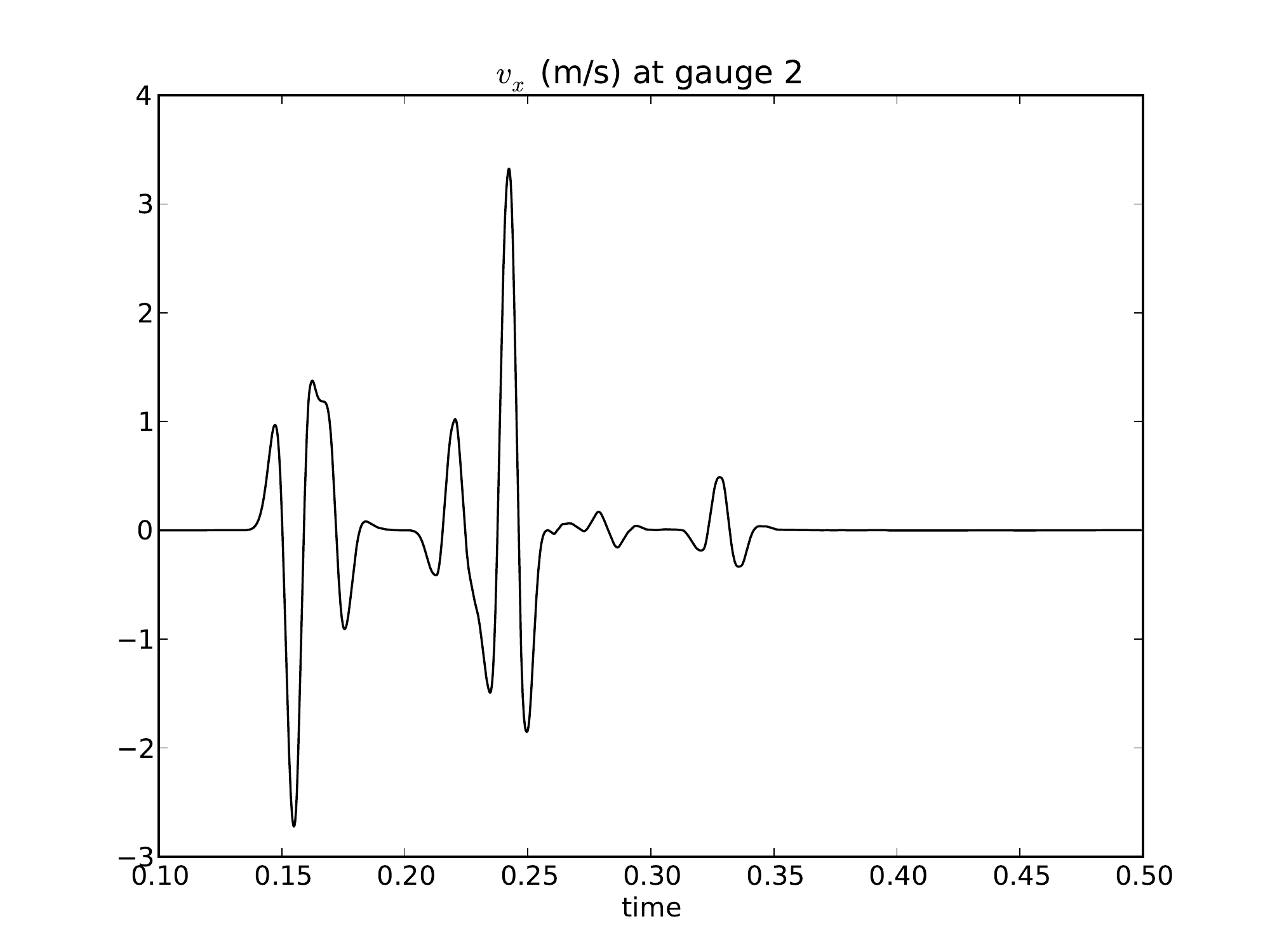}}
    \subfloat[]{\includegraphics[width=0.3\textwidth,
        trim=0.34in 0 0.6in 0, clip]{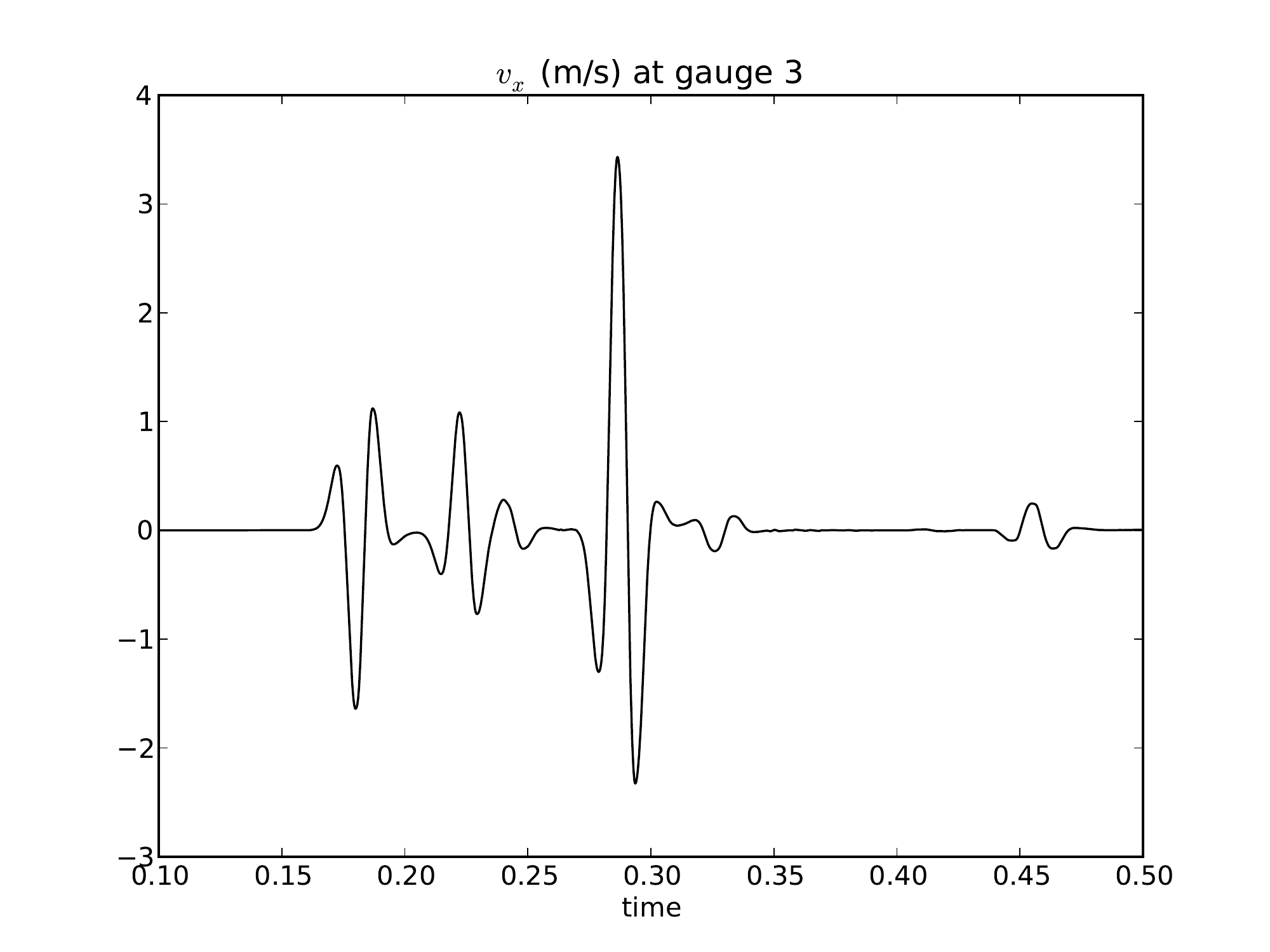}}\\
    \subfloat[]{\includegraphics[width=0.3\textwidth, page=17,
        trim=0.8in 7in 4.28in 1.2in, clip]{delapuente-2008-poro-dg-modified.pdf}}
    \subfloat[]{\includegraphics[width=0.3\textwidth, page=17,
        trim=0.8in 4.19in 4.28in 4.01in, clip]{delapuente-2008-poro-dg-modified.pdf}}
    \subfloat[]{\includegraphics[width=0.3\textwidth, page=17,
        trim=0.8in 1.38in 4.28in 6.82in, clip]{delapuente-2008-poro-dg-modified.pdf}}\\
    \subfloat[]{\includegraphics[width=0.3\textwidth,
        trim=0.34in 0 0.6in 0, clip]{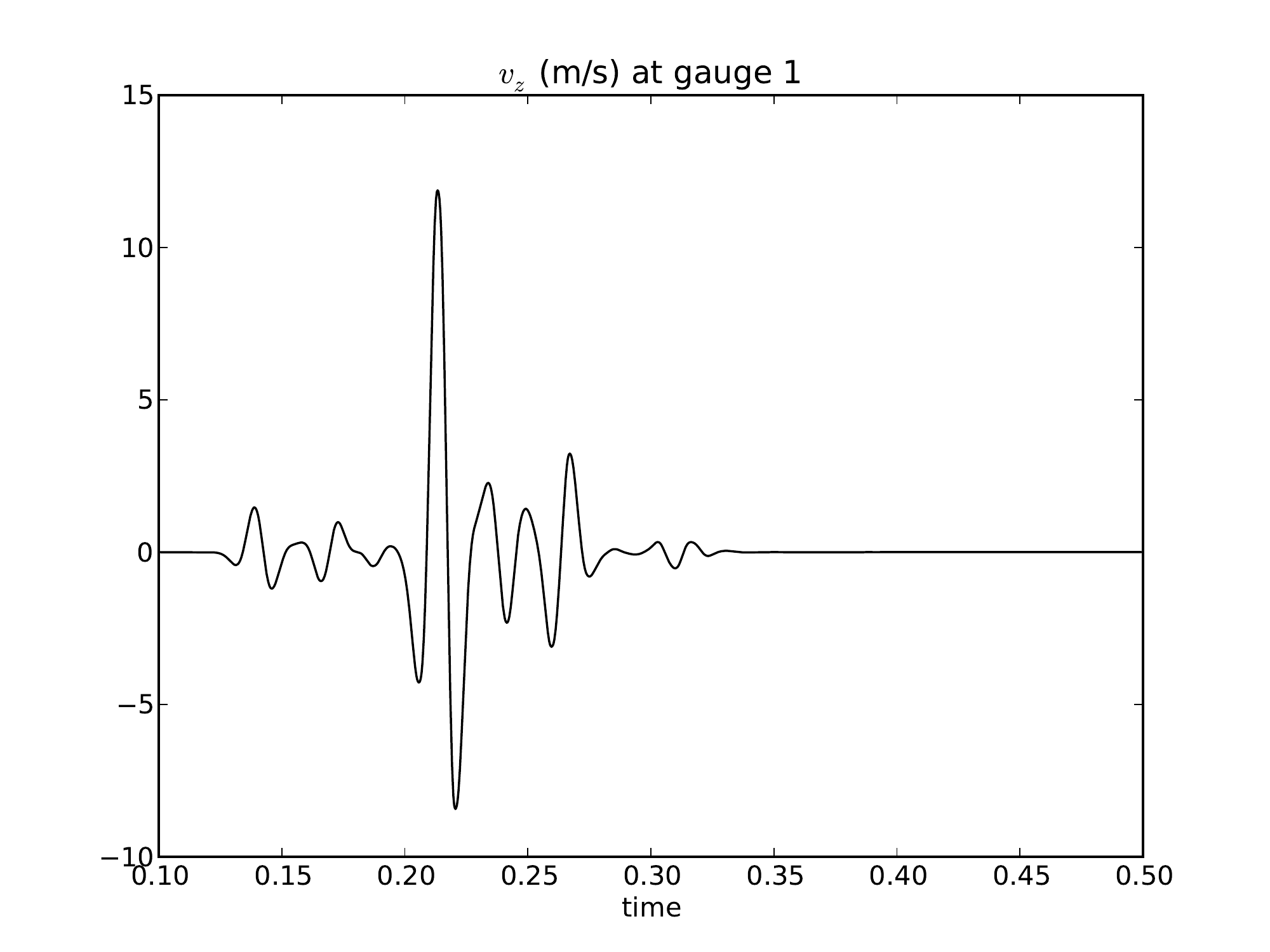}}
    \subfloat[]{\includegraphics[width=0.3\textwidth,
        trim=0.34in 0 0.6in 0, clip]{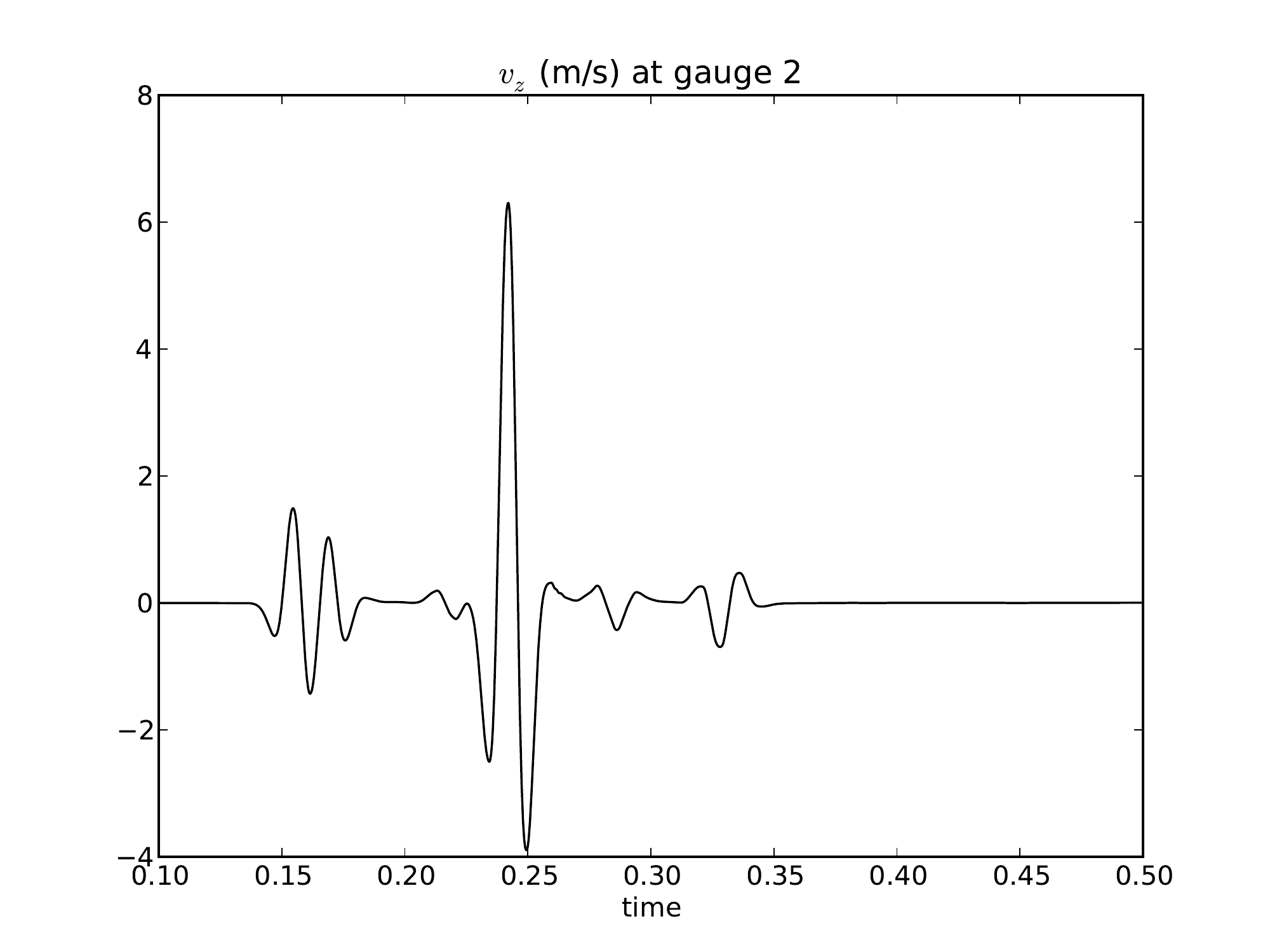}}
    \subfloat[]{\includegraphics[width=0.3\textwidth,
        trim=0.34in 0 0.6in 0, clip]{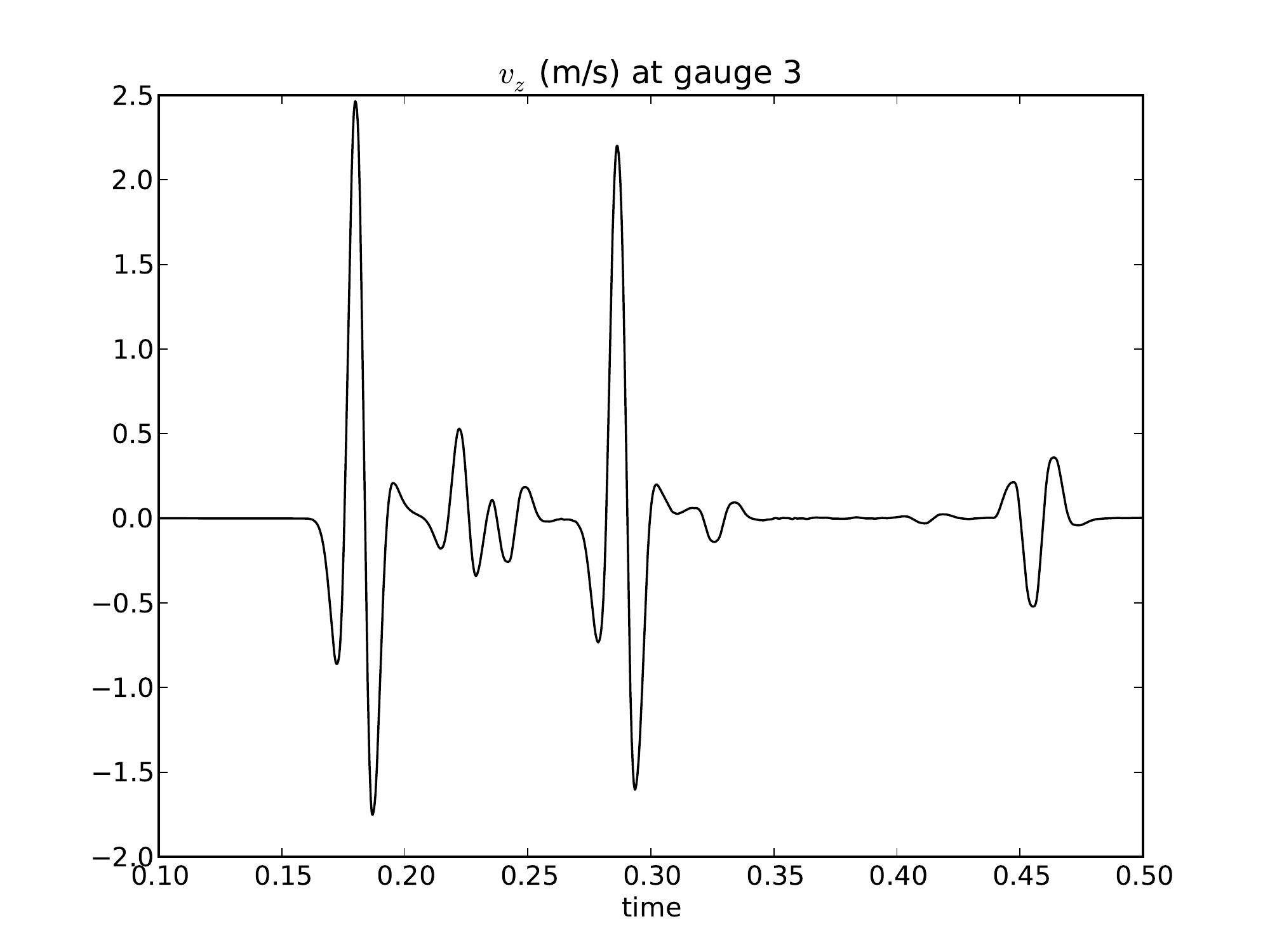}}\\
    \subfloat[]{\includegraphics[width=0.3\textwidth, page=17,
        trim=4.27in 7in 0.81in 1.2in, clip]{delapuente-2008-poro-dg-modified.pdf}}
    \subfloat[]{\includegraphics[width=0.3\textwidth, page=17,
        trim=4.27in 4.19in 0.81in 4.01in, clip]{delapuente-2008-poro-dg-modified.pdf}}
    \subfloat[]{\includegraphics[width=0.3\textwidth, page=17,
        trim=4.27in 1.38in 0.81in 6.82in, clip]{delapuente-2008-poro-dg-modified.pdf}}
    \caption{Time-histories of solid velocity at the gauges indicated
      in Figure \ref{fig:2mat-snapshot-vz}.  Top to bottom:
      $x$-direction from \clawpack{}; $x$-direction from de la Puente
      et al.~\cite{delapuente-dumbser-kaser-igel:poro-dg};
      $z$-direction from \clawpack{}; $z$-direction from de la Puente
      et al.  Left to right: gauge 1, at $(x,z) = (950\,\text{m},
      750\,\text{m})$; gauge 2, at $(950\,\text{m}, 650\,\text{m})$;
      gauge 3, at $(950\,\text{m}, 500\,\text{m})$.  Results of de la
      Puente et al.\ are reproduced in accordance with the
      policies of the publishing journal; the ADER-DG and FD curves in
      these plots are discontinuous Galerkin and finite difference
      results, respectively.  \label{fig:2mat-gauges}}
  \end{center}
\end{figure}

Figure \ref{fig:2mat-snapshot} shows a snapshot of the $z$-direction
solid velocity $0.25\,\text{s}$ after the start of the simulation,
analogous to Figure 9(a) of de la Puente et
al.~\cite{delapuente-dumbser-kaser-igel:poro-dg}, along with the AMR
grids at this time.  The solid black dot indicates the source
location, and the white-centered black dots indicate the gauge
locations.  Because the eigenvectors associated with each wave family
are different in the two materials, when a wave impinges on the
material interface it produces reflected and transmitted waves in each
of the three families; this results in a rich and complex solution
structure.  In addition to the reflected and transmitted waves, we
also see head waves, where a wave in the lower half of the domain
excites a slower wave family in the upper half.  This results in a
straight wavefront, rather than a curved one.  At this point in the
simulation the level 2 AMR grids have expanded to cover most of the
domain, but the level 3 grids are concentrated around the wavefronts.
Figure \ref{fig:2mat-gauges} shows the time histories of the solid $x$
and $z$ velocities at the three gauges.  The results are in generally
good agreement with Figure 10 of de la Puente et
al.~\cite{delapuente-dumbser-kaser-igel:poro-dg}, also shown in Figure
\ref{fig:2mat-gauges}, although the peaks
of slow P wave event at $t = 0.45 \, \text{s}$ at gauge 3 are clipped
in our simulation because of the limiter, and the magnitude of the
second large excursion in vertical velocity at gauge 3 in our solution
seems somewhat less.

While the variety of different wave speeds and reflected/transmitted
waves in this problem make adaptive mesh refinement less useful than
is typically the case --- by halfway through the simulation time,
Figure \ref{fig:2mat-snapshot-grids} shows a large portion of the
domain refined at the finest level, because there are wavefronts
present throughout the domain --- we still realize a substantial
savings in computation time.  On an Amazon EC2 Cluster 8XL instance,
running with 32 OpenMP threads, these results took 20 minutes 31
seconds to obtain, whereas a uniformly refined grid with the same cell
size as the finest AMR grids took 47 minutes 19 seconds and produced
no significant change in the solution.  (Both times are the average of
two runs; each pair differed by 4 seconds or less.)  The large number
of hardware threads available on this type of EC2 instance is the reason why
there are so many separate fine grids in Figure
\ref{fig:2mat-snapshot-grids} --- \amrclaw{} uses a coarse-grained
parallelization strategy, with each grid at each time step processed
by a single thread, which means that many grids must be present in
order to take full advantage of highly parallel computers.  The number
of grids used at each refinement level is indirectly
controlled by setting the maximum size of the individual grids; for
the AMR computation shown above, grids were allowed to extend no
more than 60 cells in any direction.  Having a very large number of
small grids (1047 level 3 grids in the figure) also eases load
balancing between threads, since most grids are of similar size, and
each thread processes many grids.

\section{Summary and Conclusions}
\label{sec:conclusion}

We have demonstrated a high-resolution finite volume code for modeling
wave propagation in porous media using Biot poroelasticity theory,
with the ability to model inhomogeneous domains and use adaptive mesh
refinement to improve solution accuracy at substantially lower
computational cost than for a uniformly refined grid.  We included the
dissipative source term in Biot's equations using operator splitting.
While this technique has produced spurious solutions when applied to
certain types of stiff source terms in the past, we have presented an
heuristic argument, based on the wave speeds of a reduced system
satisfying a certain subcharacteristic condition, that we should not
expect to encounter spurious solutions here.  This expectation was
borne out by our numerical results.

For the inviscid and viscous high-frequency regimes, our numerical
solutions converge to analytic plane wave solutions for fast P and S
waves with second-order accuracy.  Convergence rates were
somewhat impaired for the inviscid slow P wave cases tested due to the
short wavelength at the frequency chosen, which caused the waves to be
underresolved on the coarsest grids.  However, the viscous,
high-frequency slow P wave test cases showed unambiguous second-order
convergence, which indicates that we can also expect second order in
the inviscid case when the slow P wave is well-resolved.  Due to the
relative stiffness of the source term compared to the problem
timescale for low-frequency waves, we obtained only roughly first-order
accuracy for fast P and S waves in the low-frequency viscous regime.
We obtained second-order accuracy for low-frequency slow P waves, but
the slow P test cases were not directly comparable to the fast P and S
cases.  The other test cases examined, involving results for a point
source in either a homogeneous orthotropic medium, or in a layered bed
of two distinct isotropic media, agreed well with results for the same
test cases published by other authors.

There are substantial opportunities for future work based on what we
have presented here.  One possibility is the extension of the
finite volume methods used here to logically rectangular mapped grids.
This extension is straightforward, and allows modeling of more complex
domains with internal boundaries, so long as a (not necessarily
smooth) mapping function can be found that maps the internal and
external boundaries to rectangles in the computational domain.
Another opportunity for future work is the implementation of a more
accurate solution procedure at low frequencies, such as one based on
the methods discussed by Hittinger \cite{hittinger:thesis} or Pember
\cite{pember:stiff-relaxation-1993}.  We intend to explore both these
avenues in subsequent publications, along with applying the software
developed here to some specific problems.

To aid in the reproducibility of the results presented here, we
provide all of the code used to generate them at \archivelink{}.

\bibliographystyle{siam}
\bibliography{poro}

\begin{thebibliography}{10}

\bibitem{lapack-users-guide}
{\sc E.~Anderson, Z.~Bai, C.~Bischof, S.~Blackford, J.~Demmel, J.~Dongarra,
  J.~Du~Croz, A.~Greenbaum, S.~Hammarling, A.~McKenney, and D.~Sorensen}, {\em
  {LAPACK} Users' Guide}, Society for Industrial and Applied Mathematics,
  Philadelphia, PA, third~ed., 1999.

\bibitem{attenborough-berry-chen:scattering}
{\sc Keith Attenborough, David~L. Berry, and Yu~Chen}, {\em Acoustic scattering
  by near-surface inhomogeneities in porous media.}, tech. report, Defense
  Technical Information Center OAI-PMH Repository [{\tt
  http://stinet.dtic.mil/oai/oai}] (United States), 1998.

\bibitem{mjb-rjl:amrclaw}
{\sc M.~J. Berger and R.~J. LeVeque}, {\em Adaptive mesh refinement using
  wave-propagation algorithms for hyperbolic systems}, SIAM Journal on
  Numerical Analysis, 35 (1998), pp.~2298--2316.

\bibitem{biot:56-1}
{\sc M.~A. Biot}, {\em Theory of propagation of elastic waves in a
  fluid-saturated porous solid. {I}. {L}ow-frequency range}, Journal of the
  Acoustical Society of America, 28 (1956), pp.~168--178.

\bibitem{biot:56-2}
\leavevmode\vrule height 2pt depth -1.6pt width 23pt, {\em Theory of
  propagation of elastic waves in a fluid-saturated porous solid. {II}.
  {H}igher frequency range}, Journal of the Acoustical Society of America, 28
  (1956), pp.~179--191.

\bibitem{biot:62}
\leavevmode\vrule height 2pt depth -1.6pt width 23pt, {\em Mechanics of
  deformation and acoustic propagation in porous media}, Journal of Applied
  Physics, 33 (1962), pp.~1482--1498.

\bibitem{buchanan-gilbert:bone-params-hifreq-1-2007}
{\sc James~L. Buchanan and Robert~P. Gilbert}, {\em Determination of the
  parameters of cancellous bone using high frequency acoustic measurements},
  Mathematical and Computer Modelling, 45 (2007), pp.~281--308.

\bibitem{buchanan-gilbert:bone-params-hifreq-2-2007}
\leavevmode\vrule height 2pt depth -1.6pt width 23pt, {\em Determination of the
  parameters of cancellous bone using high frequency acoustic measurements
  {II}: inverse problems}, Journal of Computational Acoustics, 15 (2007),
  pp.~199--220.

\bibitem{buchanan-gilbert-khashanah:bone-params-lofreq-2004}
{\sc James~L. Buchanan, Robert~P. Gilbert, and Khaldoun Khashanah}, {\em
  Determination of the parameters of cancellous bone using low frequency
  acoustic measurements}, Journal of Computational Acoustics, 12 (2004),
  pp.~99--126.

\bibitem{bgwx:2004}
{\sc J.~L. Buchanan, R.~P. Gilbert, A.~Wirgin, and Y.~S. Xu}, {\em Marine
  acoustics: direct and inverse problems}, SIAM, Philadelphia, 2004.

\bibitem{cada-torrilhon:3rdorder-ldlr}
{\sc M.~{\v{C}}ada and M.~Torrilhon}, {\em Compact third-order limiter
  functions for finite volume methods}, Journal of Computational Physics, 228
  (2009), pp.~4118--4145.

\bibitem{carcione:poro-aniso-1996}
{\sc J.~M. Carcione}, {\em Wave propagation in anisotropic, saturated porous
  media: plane-wave theory and numerical simulation}, Journal of the Acoustical
  Society of America, 99 (1996), pp.~2655--2666.

\bibitem{carcione:wave-book}
{\sc J.~M. Carcione}, {\em Wave Fields in Real Media: Wave Propagation in
  Anisotropic, Anelastic, and Porous Media}, Elsevier, Oxford, 2001.

\bibitem{carcione-morency-santos:comp-poro-rev}
{\sc J.~M. Carcione, C.~Morency, and J.~E. Santos}, {\em Computational
  poroelasticity -- a review}, Geophysics, 75 (2010), pp.~75A229--75A243.

\bibitem{carcione-quirogagoode:biot-fdm-opsplit}
{\sc J.~M. Carcione and G.~Quiroga-Goode}, {\em Some aspects of the physics and
  numerical modeling of {B}iot compressional waves}, Journal of Computational
  Acoustics, 3 (1995), pp.~261--280.

\bibitem{chen-levermore-liu:stiff-relaxation}
{\sc G.-Q. Chen, C.~D. Levermore, and T.-P. Liu}, {\em Hyperbolic conservation
  laws with stiff relaxation terms and entropy}, Communications in Pure and
  Applied Mathematics, 47 (1994), pp.~787--830.

\bibitem{chiavassa-lombard:poro-fdm-2011}
{\sc G.~Chiavassa and B.~Lombard}, {\em Time domain numerical modeling of wave
  propagation in {2D} heterogeneous porous media}, Journal of Computational
  Physics, 230 (2011), pp.~5288--5309.

\bibitem{colella-majda-roytburd:reacting-shock}
{\sc P.~Colella, A.~Majda, and V.~Roytburd}, {\em Theoretical and numerical
  structure for reacting shock waves}, SIAM Journal on Scientific and
  Statistical Computing, 7 (1986), pp.~1059--1080.

\bibitem{cowin:bone-poro}
{\sc S.~C. Cowin}, {\em Bone poroelasticity}, Journal of Biomechanics, 32
  (1999), pp.~217--238.

\bibitem{cowin-cardoso:fabric-2010}
{\sc S.~C. Cowin and L.~Cardoso}, {\em Fabric dependence of bone ultrasound},
  Acta of Bioengineering and Biomechanics, 12 (2010).

\bibitem{cowin-mehrabadi:elastic-symmetry}
{\sc S.~C. Cowin and M.~M. Mehrabadi}, {\em Identification of the elastic
  symmetry of bone and other materials}, Journal of Biomechanics, 22 (1989),
  pp.~503--515.

\bibitem{dai-vafidis-kanasewich:poro-vel-stress}
{\sc N.~Dai, A.~Vafidis, and E.~Kanasewich}, {\em Wave propagation in
  heterogeneous porous media: a velocity-stress, finite-difference method},
  Geophysics, 60 (1995), pp.~327--340.

\bibitem{delapuente-dumbser-kaser-igel:poro-dg}
{\sc J.~de~la Puente, M.~Dumbser, M.~K{\"a}ser, and H.~Igel}, {\em
  Discontinuous {G}alerkin methods for wave propagation in poroelastic media},
  Geophysics, 73 (2008), pp.~T77--T97.

\bibitem{degrande-deroeck:poro-spectral-freq}
{\sc G.~Degrande and G.~De~Roeck}, {\em {FFT}-based spectral analysis
  methodology for one-dimensional wave propagation in poroelastic media},
  Transport in Porous Media, 9 (1992), pp.~85--97.

\bibitem{deresiewicz-skalak:poro-uniqueness-1963}
{\sc H.~Deresiewicz and R.~Skalak}, {\em On uniqueness in dynamic
  poroelasticity}, Bulletin of the Seismological Society of America, 53 (1963),
  pp.~783--788.

\bibitem{detournay-cheng:borehole}
{\sc E.~Detournay and A.~H-D. Cheng}, {\em Poroelastic response of a borehole
  in a non-hydrostatic stress field}, International Journal of Rock Mechanics
  and Mining Sciences and Geomechanics Abstracts, 25 (1988), pp.~171--182.

\bibitem{garg-nayfeh-good:poro-1974}
{\sc S.~K. Garg, A.~H. Nayfeh, and A.~J. Good}, {\em Compressional waves in
  fluid-saturated elastic porous media}, Journal of Applied Physics, 45 (1974),
  pp.~1968--1974.

\bibitem{gibson:composites}
{\sc R.~F. Gibson}, {\em Principles of Composite Material Mechanics},
  McGraw-Hill, New York, 1994.

\bibitem{gilbert-guyenne-ou:bone-2012}
{\sc R.~P. Gilbert, P.~Guyenne, and M.~Yvonne Ou}, {\em A quantitative
  ultrasound model of the bone with blood as the interstitial fluid},
  Mathematical and Computer Modelling, 55 (2012), pp.~2029--2039.

\bibitem{gilbert-lin:1997}
{\sc R.~P. Gilbert and Z.~Lin}, {\em Acoustic field in a shallow, stratified
  ocean with a poro-elastic seabed}, Zeitschrift f{\"u}r Angewandte
  {M}athematik und {M}echanik, 77 (1997), pp.~677--688.

\bibitem{gilbert-ou:2003}
{\sc R.~P. Gilbert and M.~Yvonne Ou}, {\em Acoustic wave propagation in a
  composite of two different poroelastic materials with a very rough periodic
  interface: a homogenization approach}, International Journal for Multiscale
  Computational Engineering, 1 (2003).

\bibitem{gurevich-schoenberg:interface-1999}
{\sc Boris Gurevich and Michael Schoenberg}, {\em Interface conditions for
  {B}iot's equations of poroelasticity}, Journal of the Acoustical Society of
  America, 105 (1999), pp.~2585--2589.

\bibitem{hassanzadeh:poro-1991}
{\sc S.~Hassanzadeh}, {\em Acoustic modeling in fluid-saturated porous media},
  Geophysics, 56 (1991), pp.~424--435.

\bibitem{hittinger:thesis}
{\sc J.~A. Hittinger}, {\em Foundations for the generalization of the {G}odunov
  method to hyperbolic systems with stiff relaxation source terms}, PhD thesis,
  University of Michigan, 2000.

\bibitem{kemm:limiters}
{\sc F.~Kemm}, {\em A comparative study of {TVD}-limiters -- well-known
  limiters and an introduction of new ones}, International Journal for
  Numerical Methods in Fluids, 67 (2011), pp.~404--440.

\bibitem{rjl:fvm-book}
{\sc R.~J. LeVeque}, {\em Finite Volume Methods for Hyperbolic Problems},
  Cambridge University Press, New York, 2002.

\bibitem{rjl-yee:stiff-split}
{\sc R.~J. LeVeque and H.~C. Yee}, {\em A study of numerical methods for
  hyperbolic conservation laws with stiff source terms}, Journal of
  Computational Physics, 86 (1990), pp.~187--210.

\bibitem{liu-tadmor:3rdorder}
{\sc X.-D. Liu and E.~Tadmor}, {\em Third order nonoscillatory central scheme
  for hyperbolic conservation laws}, Numerische Mathematik, 79 (1998),
  pp.~397--425.

\bibitem{mikhailenko:poro-1985}
{\sc B.~G. Mikhailenko}, {\em Numerical experiment in seismic investigations},
  Journal of Geophysics, 58 (1985), pp.~101--124.

\bibitem{morency-tromp:poro-spectral-time}
{\sc C.~Morency and J.~Tromp}, {\em Spectral-element simulations of wave
  propagation in porous media}, Geophysical Journal International, 179 (2008),
  pp.~1148--1168.

\bibitem{pember:spurious-solns}
{\sc R.~B. Pember}, {\em Numerical methods for hyperbolic conservation laws
  with stiff relaxation {I}. {S}purious solutions}, SIAM Journal on Applied
  Mathematics, 53 (1993), pp.~1293--1330.

\bibitem{pember:stiff-relaxation-1993}
\leavevmode\vrule height 2pt depth -1.6pt width 23pt, {\em Numerical methods
  for hyperbolic conservation laws with stiff relaxation {II}. {H}igher-order
  {G}odunov methods}, SIAM Journal on Scientific Computing, 14 (1993),
  pp.~824--859.

\bibitem{santos-orena:poro-fem}
{\sc J.~E. Santos and E.~J. Ore\~na}, {\em Elastic wave propagation in
  fluid-saturate porous media, part {II}: {T}he {G}alerkin procedures},
  Mathematical Modeling and Numerical Analysis, 20 (1986), pp.~129--139.

\bibitem{sharma:dissimilar-boundary-2008}
{\sc M.~D. Sharma}, {\em Wave propagation across the boundary between two
  dissimilar poroelastic solids}, Journal of Sound and Vibration, 314 (2008),
  pp.~657--671.

\bibitem{claw.org}
{\sc {The {\sc clawpack} authors}}, {\em {\sc clawpack} software}.
\newblock {\tt www.clawpack.org}.

\end{thebibliography}

\end{document}